\documentclass[11pt]{article}

\usepackage{indentfirst}
\usepackage{bbm}
\usepackage{mathtools}
\usepackage{xcolor,soul}
\usepackage{ulem}

\usepackage{amsmath}
\usepackage{amssymb}
\usepackage{latexsym}
\usepackage{amsmath, amsfonts,amssymb, amsthm, euscript,makeidx,color,mathrsfs}
\usepackage{indentfirst}
\newtheorem{assumption}{Assumption}
\def\qed{ \ \vrule width.2cm height.2cm depth0cm\smallskip}

\newcommand{\ba}{\begin{array}}
\newcommand{\ea}{\end{array}}
\newcommand{\be}{\begin{equation}}
\newcommand{\ee}{\end{equation}}
\newcommand{\bea}{\begin{eqnarray}}
\newcommand{\eea}{\end{eqnarray}}
\newcommand{\beaa}{\begin{eqnarray*}}
\newcommand{\eeaa}{\end{eqnarray*}}

\def\neg{\negthinspace}

\def\a{\alpha}

\def\e{\varepsilon}

\def\l{\lambda}
\def\m{\mu}
\def\n{\nu}

\def\si{\sigma}
\def\t{\tau}
\def\f{\varphi}
\def\th{\theta}
\def\o{\omega}

\def\f{\phi}

\def\D{\Delta}
\def\G{\Gamma}

\def\O{\Omega}

\def\G{\Gamma}
\def\D{\Delta}

\def\O{\Omega}
\def\cA{{\cal A}}
\def\cB{{\cal B}}

\def\cF{{\cal F}}
\def\cG{{\cal G}}

\def\hC{\mathbb{C}}

\def\hE{\mathbb{E}}
\def\hF{\mathbb{F}}
\def\hG{\mathbb{G}}

\def\hL{\mathbb{L}}

\def\hN{\mathbb{N}}

\def\hP{\mathbb{P}}
\def\hQ{\mathbb{Q}}
\def\hR{\mathbb{R}}

\def\hX{\mathbb{X}}

\def\sA{\mathscr{A}}
\def\sB{\mathscr{B}}

\def\sE{\mathscr{E}}
\def\scF{\mathscr{F}}

\def\scL{\mathscr{L}}

\def\sP{\mathscr{P}}

\def\sU{\mathscr{U}}

\def\no{\noindent}

\def\ss{\smallskip}
\def\ms{\medskip}
\def\bs{\bigskip}
\def\q{\quad}
\def\qq{\qquad}

\def\pa{\partial}
\def\cd{\cdot}
\def\cds{\cdots}
\def\lan{\langle}
\def\ran{\rangle}

\def\bx{{\bf x}}

\def\tr{\hbox{\rm tr}}

\def\qed{ \hfill \vrule width.25cm height.25cm depth0cm\smallskip}
\newcommand{\dfnn}{\stackrel{\triangle}{=}}
\newcommand{\basa}{\begin{assumption}}
\newcommand{\easa}{\end{assumption}}

\newcommand{\bas}{\begin{assum}}
\newcommand{\eas}{\end{assum}}

\def\lan{\mathop{\langle}}
\def\ran{\mathop{\rangle}}

\def\pa{\partial}

 \def\cd{\cdot}
\def\cds{\cdots}

\def\as{\hbox{\rm-a.s.{ }}}

\def\tr{\hbox{\rm tr$\,$}}

\def\ind{{\perp \neg\neg\neg\perp}}

\def\dis{\displaystyle}

\def\bx{{\bf x}}
\def\bv{{\bf v}}

\def\1{{\bf 1}}

\def\:{\!:\!}

at 9pt

\def\prehp(#1,#2){\ensuremath{  #1 \cdot #2 }}

\begin{document}

\newtheorem{thm}{Theorem}[section]
\newtheorem{lem}[thm]{Lemma}
\newtheorem{cor}[thm]{Corollary}
\newtheorem{prop}[thm]{Proposition}
\newtheorem{rem}[thm]{Remark}
\newtheorem{eg}[thm]{Example}
\newtheorem{defn}[thm]{Definition}
\newtheorem{assum}[thm]{Assumption}

\renewcommand {\theequation}{\arabic{section}.\arabic{equation}}
\def\thesection{\arabic{section}}

\title{\bf A Generalized Kyle-Back Strategic Insider Trading Model with Dynamic Information}

\author{
Jin Ma\thanks{ \noindent Department of
Mathematics, University of Southern California, Los Angeles, 90089;
email: jinma@usc.edu. This author is supported in part
by US NSF grants  \#1106853. } 
~ and ~ Ying Tan\thanks{\noindent Department of
Mathematics, University of Southern California, Los Angeles, 90089; email: yingtan@usc.edu}}

\date{}
%{\today}
\maketitle

\begin{abstract}
In this paper we consider a class of generalized Kyle-Back strategic insider trading models in which the insider is able to use the dynamic information obtained by observing the instantaneous movement of an underlying 
asset that is allowed to 
be influenced by its market price. Since such a model will be largely outside the Gaussian paradigm, we shall try to Markovize it by introducing an auxiliary diffusion process, in the spirit of the {\it weighted total order process} of, e.g., \cite{CCD11}, as a part of the ``pricing rule". As the main technical tool in solving the Kyle-Back equilibrium, we study a class of  {\it Stochastic Two-Point Boundary Value Problem} (STPBVP), which 
resembles the {\it dynamic Markov bridge} in the literature, but without insisting on its local martingale requirement. 
In the case when the solution of the STPBVP has an {\it affine structure},
we show that the pricing rule functions, whence the Kyle-Back equilibrium, can be determined by the decoupling field of a {\it forward-backward SDE} obtained via a non-linear filtering approach, along with a set of compatibility conditions.
\end{abstract}

\vfill \bs

\no{\bf Keywords.} \rm Kyle-Back equilibrium, strategic insider trading problem, conditioned SDE, two-point boundary value problem, FKK equation, forward-backward SDE, stochastic optimal control.

\bs

\no{\it 2020 AMS Mathematics subject classification:} 60H10; 93E11; 91G15, 80.

\eject

\section{Introduction}
\setcounter{equation}{0}

In this paper we are interested in an asset pricing problem with asymmetric information, known as the {\it Kyle-Back strategic insider trading equilibrium problem} initiated by  Kyle \cite{K85} and  Back (\cite{KB92, BP98}) (see also  \cite{ABO, BHMO12, CalSt, Dani, HS} and the references therein for various generalizations of such models, along with different approaches).  
In particular, we will focus on the cases of {\it dynamic information}, in which the insider is allowed to use the dynamically observed information on the underlying asset, rather than the information at a fixed terminal time, as it was originally suggested. We shall carry out the analysis in a general  Markovian, hence non-Gaussian framework.

The Kyle-Back strategic insider trading problem can be briefly described as follows. Consider a market that involves three types of agents:
(i) {\it The insider}, who possesses some information of a given asset $V=\{V_t\}_{t\in[0,T]}$ that is not observable in the market. The information can be either the law of $V_T$, or the instantaneous observation of the state $V_t$, $t\in [0,T]$, or both. In the literature, they  are often referred to as the ``long-lived information" and the ``dynamic information", respectively.  
The insider will then submit her order, denoted by $\xi_t$,  $t\in[0,T]$.  (ii) {\it The noise traders}, who have no direct information of the  asset $V$, and (collectively) submit an order $z_t$ at time $t\in [0,T]$. It is commonly assumed, by virtue of the central limit theorem, 
that $z_t=\int_0^t\sigma^z_{t}dB^z_{t}$,  where $B^z$ is a Brownian motion. (iii) Finally,  the {\it marked maker},  who observes  the total traded volume in the market, $Y_t\dfnn\xi_t+z_t$, $t\in[0,T]$, and sets the price for $V_t$. It is standard to assume  (see, e.g., \cite{K85}, by a Bertrand competition argument) that the market price $P_t$, $t\ge 0$, is the $L^2$-projection of the true value $V$ to the space of $\hF^Y$-measurable random variables. In other words, one assumes that, for $t\in[0,T]$,
\bea
\label{P0}
P_t=\left\{\ba{lll}
   \hE[V_T|\cF^Y_t]  \qq &\mbox{\rm (long-lived information)}\ms\\
   \hE[V_t|\cF^Y_t]  &\mbox{(dynamic information)},
   \ea\right. 
\eea
where ${\cal F}_t^Y\dfnn \si\{Y_s, s\leq t\}$. An {\it equilibrium} of the Kyle-Back problem consists of an insider's strategy $\xi^*$ that maximizes her expected wealth at the terminal time $T$, together with a specified market price $P$ in either form of (\ref{P0}) (often referred to as the {\it market efficiency}).

Strong efforts have been made in recent years to extend the Kyle-Back problem to more general settings beyond the traditional Gaussian framework, and some deeper mathematical tools have been introduced to deal with the solvability issues accompanied by the generality of the modeling (see, for example, \cite{CCD11, CCD13a,CCD13b, CEL} and the references cited therein). It is thus always interesting to identify methodologies that are easily accessible and at the same time efficient  for solving more general models. This paper is an effort in  this general direction. 

We are interested in a Kyle-Back equilibrium problem with the following features:
 
 (i) The evolution of the dynamics of the underlying asset can depend on the market price $P=\{P_t\}$ (hence not independent of the market information $\hF^Y=\{\cF^Y_t\}$).     
 
(ii) The insider can observe both the movements of the underlying asset and the market price, and uses the information when decides  her optimal strategy; and 

(iii) the market maker's pricing rule is of the form an ``{\it optional projection}" of the underlying asset (i.e. the second form in (\ref{P0})), rather than a martingale (the first form in (\ref{P0})).

We note that the feature (i) above, although reasonable (see, e.g., \cite{MSZ}), would put our problem outside most of the cases studied in the literature, due to various technical reasons which will become clear when our analysis proceeds, especially when the idea of ``dynamic Markov bridge" is adopted. The requirement (iii), however, will be a natural connecting point to the nonlinear filtering, given the reasonable structure of the asymmetric information. More precisely, in this paper we shall assume that the underlying asset $V$ is governed by the following general SDE:
\bea
\label{V0}
dV_t= b(t, V_{\cd\wedge t}, P_{\cd\wedge t})dt +\si(t, V_{\cd\wedge t}, P_{\cd\wedge t})dB^1_t, \qq V_0=v,
\eea
where $b,\si$ are given measurable functions. We shall also assume, as commonly seen in the literature, that the insider's strategy is of the form $\xi_{t}=\int_0^t\a_sds$, $t\ge 0$, where the ``rate" $\a$ can depend on both $V$ and $P$ in an {\it nonanticipative} way, so  that  the dynamics the market maker observes is:
\vspace{-1mm}
\begin{equation}
\label{y1}
dY_t=d\xi_t+dz_t=\a(t,V_{\cd\wedge t}, P_{\cd\wedge t})dt+dB^2_t, \qq t\ge 0.
\end{equation}

We remark that under the market efficiency requirement (\ref{P0}), the SDEs (\ref{V0}) and (\ref{y1}) in general form a so-called 
{\it conditional mean-field SDE} (CMFSDE) (or more generally, {\it conditional McKean-Vlasov SDE} (CMVSDE), whose well-posedness is not trivial (cf., e.g, \cite{BLM, MSZ}). 
In this paper we shall take a different route, and follow the idea of \cite{CCD11} and introduce a {\it factor} model which in a sense {\it Markovizes} the ``path-dependent" SDEs (\ref{V0}) and (\ref{y1}) completely. To be more precise, we are looking for a {\it factor} process $X$ that is determined completely by the observation 
It\^o process $Y$, in the sense that $X_t=\Psi(t, Y_{\cd\wedge t})$, such that the market price $P$ is determined by 

\ss
\centerline{
$P_t=H(t, X_t)=H(t, \Psi(t, Y_{\cd\wedge t}))=\Phi(t, Y_{\cd\wedge t})$,  $ t\in [0,T]$.
}

\ss
\no Such a factor process $X$ resembles the so-called {\it weighted total} process (see, e.g., \cite{CCD11}), which was often assumed to be a diffusion process driven by the observation process $Y$ (see \S2 for detailed discussion). 
With such a Markovization, we shall recast the equilibrium problem as a stochastic control problem and show that, by a  dynamic programming argument, a necessary condition for $\a^*=\tilde{u}^*(\cds)$ being optimal is that the corresponding solution $(V, X)$ satisfies:
\vspace{-1mm}
\bea
\label{terminal1}
V_T=P_T=H(T, X_T):=g(X_T).
\eea
We note that the relationship (\ref{terminal1}) naturally leads to a {\it two-point boundary value problem} structure, or a ``bridge". In fact, there has been a tremendous effort to use the notion of  {\it dynamic Markov bridge}  to help finding the Kyle-Back equilibrium (see, e.g., 
\cite{FWY, CCD11, CCD13a}), and the methodology works well when some technical and structural assumptions are made to ensure the solvability. However, these assumptions excludes the more convoluted situations such as (\ref{V0}).

The main motivation of this paper is based on the following observation: although dynamic Markov bridge is a powerful tool in solving the problem, it can be slightly relaxed for the purpose for this particular problem. In other words, a slightly generalized version, which we shall refer to as the stochastic  {\it two-point boundary value problem} (STPBVP), would be sufficient, if not more effective, for our purpose. Our main idea is to simply use the so-called ``conditioned" SDE (see, Baudoin \cite{FB02}) and design a specific {\it minimal} probability measure for the two-dimensional Markovian process $(V, X)$, and construct a weak solution to the STPBVP. Some fundamental tools in the study of dynamic Markov bridge should be sufficient for the resolution of TPBVP, whence the desired Kyle-Back equilibrium problem. 
We should note that the choice of the coefficients of the factor process $X$ is somewhat {\it ad hoc}, and we can and will impose some structural assumptions that would lead to  explicit ``compatibility conditions" among coefficients of $V$ and $X$. In particular, in this paper we shall assume an  {\it affine structure}, motivated in part by the well-known Widder's Theorem (cf. e.g., \cite{BPSV, G62, W44, W53})
and the solution of the STPBVP. We shall first argue that, given the affine structure,  some analysis similar to   {\it affine term structure of interest rates} can be used to derive the compatibility conditions;
and  the optional projection $P_t=\hE[V_t|\cF^Y_t]$ can be rigorously put into a nonlinear filtering framework with $(V,X)$ being the state signal process, and $Y$ being the observation process. Furthermore, the terminal condition (\ref{terminal1}) will lead to a coupled
{\it Forward-backward SDE} (FBSDE), with the factor process $X$ being the forward SDE, and the Fujisaki-Kallianpur-Kunita (FKK) equation of the filtering problem being the backward SDE, both
driven by the process $Y$. We then show that the corresponding {\it decoupling field} (cf. \cite{MWZZ}) is exactly the pricing rule $H$ (see, e.g., \cite{CCD11}). Note that such a connection opens the door to a potentially much more general framework in which   the decoupling field $H$ is  allowed to be  a random field, determined by a {\it backward stochastic PDE} (BSPDE),  
as is often seen in the FBSDE literature (cf. e.g., \cite{MYbook}). We hope to be able to address such issues in our future publications.

The rest of the paper is organized a follows. In \S2 we formulate the problem and introduce the notations and definitions. In \S3 we revisit the conditioned SDE; and in \S4 we formulate the stochastic two-point boundary value problem (STPBVP) and investigate its well-posedness and  fundamental properties. In \S5 we introduce the notion of  affine structure for the solution to the STPBVP and associated insider strategies. In \S6 we discuss the filtering problem and derive the FKK equation and  the corresponding FBSDE under the affine structure. Finally, in \S7 we discuss the sufficient conditions for optimality, and determine the equilibrium strategies.

\section{Preliminaries and Problem Formulation}
\setcounter{equation}{0}

Throughout this paper, let $\hX$ be a  generic Euclidean space  and regardless of its dimension, $\lan\cd,\cd\ran$ and $|\cd|$  be its inner product and norm, respectively. We denote the space of $\hX$-valued continuous functions defined on $[0,T]$ with the usual sup-norm by $\hC([0,T];\hX)$. In particular, we denote $\hC^2_T:=\hC([0,T];\hR^2)$,
and  let $\sB(\hC_T^2)$ be its topological Borel field.
 We shall assume that all randomness in this paper is characterized by a {\it canonical probabilistic set-up}:
 $(\O,\cF, \hP, \hF, B)$, where $(\O, \cF):=(\hC^2_T, \sB(\hC^2_T)$; $\hP\in\sP(\O)$;   
  and  $B=(B^1, B^2)$ is a $\hP$-Brownian motion. Moreover, we shall 
assume that $\hF^i=\{\cF^{B^i}_t\}_{t\ge 0}$, $i=1,2$, is the natural filtration generated by
$B^1$ and $B^2$, respectively, and  $\hF=\hF^1\vee\hF^2$, with the usual $\hP$-augmentation so that it satisfies the {\it
usual hypotheses} (cf. e.g., \cite{PR}). Finally, we denote $\hQ^0\in \sP(\O)$ to be
 the Wiener measure on $(\O, \cF)$;  $B^0_t(\o)=\o(t)$, $\o\in\O$, the canonical process; and  $\hF^0:=\{\cF^0_t\}_{t\in[0,T]}$,  where $\cF^0_t:=\sB_t(\hC_T^2):=\si\{\o(\cd\wedge t):\o\in\hC^2_T\}$, $t\in[0,T]$.
In what follows we shall make use of the following notations:

\ss
$\bullet$ For any sub-$\si$-field $\cG\subseteq\cF_T$ and $1\le p<\infty$,
$L^p(\cG;\hX)$ denotes the space of all $\hX$-valued, $\cG$-measurable
random variables $\xi$ such that $\hE|\xi|^p<\infty$. As usual, $\xi\in L^\infty
(\cG;\hX)$ means that it is $\cG$-measurable and bounded.

\ss

$\bullet$  For $1\le p<\infty$, $\hG\subseteq \hF$, $L^p_\hG([0,T];\hX)$ denotes the space of all
$\hX$-valued, $\hG$-progressively measurable processes $\xi$ satisfying
$\hE\int_0^T|\xi_t|^pdt<\infty$. The meaning of $L^\infty _\hG([0,T];\hX)$
is defined similarly.
For simplicity, we will often drop $\hX(=\hR)$ from the notation, and denote all ``$L^p$-norms" by
$\|\cd \|_p$, regardless it is for $L^p(\cG)$, or for $L^p_{\hF}([0,T])$, when the context is clear. 

\ms
{\bf The Problem Formulation.} As we indicated in before, there are three types of agents in the market: the insider; the noise trader; and the market maker, which we now  specify in details.

\ss
(i) {\it The insider}. 
In this paper we shall assume that the insider can both dynamically observe the liquidation value of the underlying asset $V=\{V_t\}$, and have some information of $V_T$, in particular, the law of $V_T$, denoted by  $m^*\in\sP(\hR)$. Specifically, we assume that the asset process $V$ is governed by 
the following SDE:
 \bea
 \label{V1}
 dV_t = b(t, V_t, P_t)dt + \sigma(t, V_t, P_t) dB^1_t, \qq & V_0=v,
 \eea
 where $ b, \si $ are measurable functions, and $P=\{P_t\}$ is the market price. We should note that  
 allowing $(b, \si)$ to depend on the market price $P$ is one of the main features of this paper, which amounts to saying that the fundamental price $V$ 
 is convoluted with the market information $\hF^Y$ (see (\ref{y}) below), which leads to some fundamental difficulties that distinguishes this paper from most of the existing literature, especially in terms of the {\it dynamic Markov bridge}.

We should note that although the insider has more information of the underlying asset, even it's law at a future time, we shall insist that its strategy is in the non-anticipating manner. More precisely, we shall assume that the order process $\xi_t$, $t\in[0,T]$, takes the form
 $\xi_t =\xi^\a_t:=\int_0^t \a_s ds$, where the process $\a=\{\a_t\}$ is called the  {\it intensity} of the trading strategy, and is   assumed to have the form $\a_t=u(t, V_{\cd\wedge t}, P_{\cd\wedge t})$, $t\in[0,T]$, for some function $u$  to be determined (see, e.g., \cite{BP98, MSZ}).

\ss

(ii) {\it The noise traders}. For simplicity, in this paper we shall assume that the (collective) order submitted by the noise traders is simply the $z_t=B^2$, for some Brownian motion $B^2\ind B^1$. In other words, we assume that $B^z=B^2$, and $\si^z\equiv 1$. 

\ss
(iii) {\it The market maker}. By virtue of the so-called Bertrand competition argument (see, e.g., \cite{K85}), we assume that at each time $t\in[0,T]$, the market maker sets the (market) price 
$P_t$ to be the ($L^2$-)projection of the (unobservable) underlying price $V_t$ onto the space of all $\cF^Y_t$-measurable random variables. That is,  
$P_t = \hE[V_t | \cF_t^Y]$, $t\in[0,T]$, where $Y$ is the total trading volume:
\begin{equation}
\label{y}
Y_t=\xi^\a_{t}+B^2_t = \int_0^t\alpha_s ds+ B^2_t, 
\qq t\in[0,T].
\end{equation}
 
 Furthermore, we require that the asymmetry of information ends at the terminal time $T$. That is, at terminal $T>0$ the value of the underlying asset $V_T$ will be revealed (by, e.g., an announcement, cf. \cite{BP98}) and  the market  price will be set as
$P_T=V_T$, so that the insider does not have any information advantage by the time $T$. We should note that such a 
requirement 
is not a natural consequence given the market parameters (i.e., the coefficients of SDEs involved), but rather one of the conditions the equilibrium  strategy must satisfy. 

Before we formulate the equilibrium problem, let us specify the set of {\it admissible strategies}:
\bea
\label{Uad}
\sU_{ad}:= \{\a\in \hL^2_{\hF}([0,T]): L^\a \mbox{~is a local martingale on
$[0,T)$}\}.
\eea
where $L^\a $ is the Dol\'eans-Dade stochastic exponential: $L^\a_t:=\exp\big\{\int_0^t\a_s dB^2_s-\frac12\int_0^t\vert\a_s\vert^2ds\big\}$, $t\in [0,T)$. 
A (generalized) Kyle-Back equilibrium consists of  a ``pricing rule" $P_t=\hE[V_t|\cF^Y_t]$, $t\in[0,T]$; and an optimal strategy $\a^* \in \sU_{ad}$, such that the terminal wealth, defined by  
\beaa
W_T= W^{\a^*}_T:=\int_0^T \xi^{\a^*}_t dP_t
\eeaa
has a maximum expected value $\hE^\hP[W^{\a^*}_T]=\sup_{\a\in\sU_{ad}}\hE^\hP[W^\a]$. 
\begin{rem}
\label{remark0} {\rm
(i) We note that, in (\ref{Uad}) the process $L^\a$ is defined only on $[0, T)$. Indeed, in light of the existing results, 
the optimal strategy $\a_t$ may very well explode when $t\nearrow T$, because the insider will try to use all the information advantage before it ends. 

\ss
 (ii) From (\ref{y}) we observe that $Y$ depends on the choice of $\a$, thus so does the process $P$, 
whence the asset price $V$. Denoting $V=V^\a$, a more precise definition of the admissible control set
 should be all $\a \in\sU_{ad}$ such that $V^\a_T\sim m^*\in\sP(\hR)$, the law that the insider was expecting. We  prefer not to impose such a restriction in order to avoid unnecessary technical subtlety, but will emphasize this issue when it is needed in our discussion (e.g., in \S7). 
\qed}
\end{rem}

{\bf The Markovization.} We note that the market price $P_t=\hE[V_t|\cF^Y_t]$, $t\in[0,T]$, is in general an {\it optional projection} of $V$ onto the filtration $\hF^Y=\{\cF_t\}$, but not necessarily an 
$\hF^Y$-martingale 
as the ``long-lived information" case (see ({\ref{P0})) considered in most of the existing literature. In general the market price $P$ can be written as 
$P_t =\Phi(t, Y_{\cd\wedge t})$,  $t\ge0$, for  some measurable function $\Phi$  defined on 
$\hC([0,T])$. Therefore (\ref{V1})--(\ref{y}) is by nature a system of ``path-dependent"
{\it Conditional McKean-Vlasov SDEs} (CMVSDEs) or
{\it Conditional Mean-field SDEs} (CMFSDEs) (see \cite{BLM, MSZ}). 
In this paper we shall follow the idea of \cite{CCD11} to first {\it Markovzie} the system (\ref{V1})-(\ref{y}) by introducing 
a {\it factor}  process $X$, which satisfies an auxiliary SDE of the form:
\bea
\label{X0}
dX_t=\m(t, X_t)dt+\rho(t, X_t)dY_t, \qq X_0=x, 
\eea
where the coefficients $(\m,\rho)$ are to be determined, so that 
the market price $P$ can be written as $P_t=H(t, X_t)$ for some function $H$. We note that, if on some probability space 
$(\O, \cF, \hQ)$, where $\hQ\in\sP(\O)$ under which $Y$ is a Brownian motion, then, as the strong solution to SDE (\ref{X0}),  $X$ can be written as $X_t=\Psi(t, Y_{\cd\wedge t})$, for some measurable function $\Psi$, and consequently, we have 
\beaa
\label{P}
P_t=\hE[V_t|\cF^Y_t]=H(t, X_t)=H(t, \Psi(t, Y_{\cd\wedge t}))=\Phi(t, Y_{\cd\wedge t}), \qq t\in [0,T].
\eeaa
We note that the factor process $X$ is similar to the {\it weighted total order process}, proposed in \cite{CCD11}), and the function $H$ (together with the coefficients $(\m, \rho)$) can be considered as the ``pricing rule" (see \cite{CCD11, CCD13a}). They will be the main 
subject of this paper.

We should remark here that a direct consequence of the Markovization is that we can now put the problem of finding the equilibrium into a standard stochastic control framework. More specifically, since $P_t = H(t,X_t)$, by a slight abuse of notation, we shall assume from now on that the underlying asset $V$ and the factor process $X$ follow  a system of  SDEs:
\bea
\label{vx00}
\left\{\ba{lll}
dV_t = b(t, V_t, X_t)dt + \sigma(t, V_t, X_t) dB^1_t,  \qq V_0=v; \ms\\
dX_t = \m(t, X_t)dt + \rho(t, X_t)dY_t=[\m(t, X_t)+\a_t\rho(t, X_t)]dt+\rho(t,X_t)dB^2_t, \qq  X_0=x.
\ea\right.
\eea
Considering (\ref{vx00}) as a controlled system with controls $\a\in\sU_{ad}$. Following the argument of \cite{KB92} by allowing  a market clearing jump at terminal time, then a simple integration by parts shows that 
the expected terminal wealth  can be written as: 
\bea 
\label{w2}
\hE[W^\a_{T}]=\hE\Big[ (V_T - P_T)\xi^\a_T + \int_0^T \xi^\a_{t}dP_t\Big]=\hE\Big[\xi^\a_TV_T-\int_0^T\a_{t}P_tdt\Big]=\hE\Big[\int_0^T[V_T-P_t]\a_{t}dt\Big]. 
\eea 
Assuming now the process $\a$ takes the feedback form: $\a_t=u(t,V_t, X_t)$, then $(V, X)$ becomes Markovian, and we deduce from (\ref{w2}) that 
\bea
\label{w1}
\hE[W^\a_{T}]=\hE\Big[\int_0^T[\hE[V_T|\cF^{V,X}_t]-P_t]\a_{t}dt\Big]=\hE\Big[ \int_0^T[ F(s,V_s, X_s) - H(s, X_s)]\a_{t}dt\Big],
\eea
where $F$ is a continuous function satisfying $F(T, v,x)=v$, and can be determined by the Kolmogorov backward equation or Feynman-Kac formula (see \S7 for details). Consequently, we can define a stochastic control problem with $(V,X)$ as the controlled dynamics, and the 
{\it cost functional}:
\bea
\label{J0}
			J(t,v,x; u) :=\hE _{t,v,x}\Big[ \int_t^T(F(s,V_s,X_s) - H(s,X_s))u(s,V_s,X_s)ds\Big],
\eea
so that the value function $\bv(t,v,x):=\sup_{\a\in \sU_{ad}}J(t,v,x; u)$ satisfies the following HJB equation:
\bea
\label{HJB0}
0&=&\pa_t \bv(t,v,x)+b(t,v,x)\pa_v\bv+ \mu(t,x)\pa_x \bv+ \dfrac{1}{2}\sigma^2(t,v,x) \pa_{vv}\bv+  \dfrac{1}{2}\rho^2(t,x) \pa_{xx}\bv\nonumber \\
&&+ \sup_{u \in \hR}\big\{   [ \rho(t,x)\pa_x \bv + F(t,v,x)-H(t,x)]u \big\}.
\eea
Clearly, a necessary condition for the ``sup"-term in (\ref{HJB0}) to be finite is 
	\beaa
		 \rho(t,x)\pa_x \bv+ F(t,v,x)-H(t,x) = 0, \qq (t, v,x) \in [0,T]\times \hR^2.
	\eeaa
 In particular, noting that $F(T,v,x) = v$, and $\bv(T, v,x)\equiv 0$ by definition (\ref{J0}), we deduce that 
    \bea
   \label{cond0}
   0\equiv \rho(T,x) \pa_x \bv(T,v,x) = H(T,x) -F(T,v,x) =:g(x)  -v, \qq(v,x) \in \hR^2,
    \eea
where  $g(x)=H(T, x)$. In other words, 
at the terminal time $T$, it holds that $V_T=g(X_T)$ for some
function $g$. In fact, similar to  \cite{CCD11}, we shall assume from now on that the function $g$ is increasing. Consequently, (\ref{cond0}) indicates an important fact:  a necessary condition for $\a\in \sU_{ad}$ being an equilibrium is that the following condition holds at the terminal time $T$:
\bea
\label{terminal}
V_T=P_T=H(T, X_T)=g(X_T).
\eea

\ss 
{\bf A Stochastic Two-Point Boundary Valued Problem  (STPBVP).} Summarizing the discussion above we see that we should look for $\a\in \sU_{ad}$ and coefficients $(\m, \rho)$ so that the following system of SDEs with initial-terminal conditions is 
solvable:
    \begin{equation}
	\label{vx}
	\left\{\begin{array}{lll}
	dV_t = b(t, V_t, X_t)dt + \sigma(t, V_t, X_t) dB^1_t, \ms\\
	dX_t = \m(t, X_t)dt + \rho(t, X_t)dY_t=[\m(t, X_t)+\a_t\rho(t, X_t)]dt+\rho(t,X_t)dB^2_t, \ms\\
      V_0=v, \q  X_0=x, \q V_T=g(X_T).
      \end{array}\right.
	\end{equation}
In what follows we shall refer to (\ref{vx}) as a {\it Stochastic Two-Point Boundary Value Problem}, which will be studied in details in the next section. The solvability of the STPBVP depends on the 
choice of the process $\a=\{\a_t\}$. In particular, we are particularly interested in the case 
when $\a$ takes the form $\a_t=u(t, V_t, X_t)$, which will render the solution $(V, X)$ a Markov process.

  We remark that the TPBVP (\ref{vx}) is closely related to the so-called {\it dynamic Markov bridge} studied in, e.g.,  \cite{FWY, CCD11, CCD13a}. In fact, if $b=\m=0$, $\si=\rho=1$, and $g(x)=x$, the problem (\ref{vx}) was first studied, as the Brownian bridge, in the
 context of insider trading in \cite{FWY}. The more general cases were considered recently in \cite{CCD11, CCD13a, CCD21}, also  in the bridge context. But on the other hand, we note that 
in the description of the problem above we see that the TPBVP (\ref{vx}) does not actually require that the solution $X$ to be a local martingale under its
own filtration, a key requirement to be a Markovian bridge (see \S3 for a more detailed discussion). Thus, the main point of this paper is 
to show that such a relaxation enables us to solve the Kyle-Back equilibrium problem in a much more general setting.

\section{The Conditioned SDE Revisited}
\setcounter{equation}{0}

Our construction of the (weak) solution to TPBVP (\ref{vx}) is based on the notion of the so-called {\it conditioned SDE} (cf. \cite{FB02}), which we now briefly describe. Recall the canonical probabilistic set-up  $(\O, \cF, \hQ^0; \hF, B^0)$ defined in the beginning of \S2. In particular, we denote the canonical process by $B^0=(B^1,Y)$ so that it is a $(\hQ^0, \hF)$-Brownian motion. Consider the SDE on canonical space $(\O, \cF, \hQ^0, B^0)$, for $t\in[0,T]$:
\bea
\label{SDEVXQ}
	\left\{\begin{array}{lll}
	dV_t = b(t, V_t, X_t)dt + \sigma(t, V_t, X_t) dB^1_t, \qq & V_0=v; \ms\\
	dX_t = \mu(t, X_t)dt +\rho(t,X_t)dY_t, & X_0=x.
	\end{array}\right.
\eea
We assume that the coefficients $b, \si, \m, \rho$ and $g$ satisfy the following  {\it Standing Assumptions}:
\begin{assum}
\label{assump1}
(i) The functions $b, \si:[0,T]\times \hR^2\mapsto \hR$ and $\m, \rho:[0,T]\times \hR \mapsto \hR$ are measurable, and continuous
in $t\in[0,T]$;

\ss
(ii) There exists $L>0$, such that, for any $t\in[0,T]$, $v, v', x,x'\in \hR$, it holds that, 
\beaa
\label{Lip}
\left\{\ba{lll}
 |b(t,0,0)|+|\si(t,0,0)|+|\m(t,0)|+|\rho(t,0)|\le L, \ms\\
|\f(t, v,x)-\f(t, v',x')|\le L(|v-v'|+|x-x'|), \qq\q & \f=b, \si, \ms\\
|\psi(t, x)-\psi(t, x')|\le L|x-x'|, &\psi=\m, \rho;
\ea\right.
\eeaa

\ss
(iii) There exists a constant $\l_0>0$, such that 
$\si(t, v,x)\ge \l_0$, $(t,v,x)\in [0,T]\times\hR^2$;

\ss
(iv) The function $g$ is uniformly Lipschitz continuous and strictly monotone increasing.
\qed
\end{assum}

Clearly, under Assumption \ref{assump1}, SDE (\ref{SDEVXQ}) has a unique strong solution over $[0,T]$,  on $(\O, \cF, \hQ^0)$, denoted by $\xi:=(V^0, X^0)$. Moreover, $\xi$ is a Markov process, and 
 we denote its transition density  by  $p(s, x; t, y)$, $0\le s<t\le T$, $x, y\in \hR^2$. Now, for any $\nu \in\sP(\hR^2)$, consider the   triplet $(T, \xi_T, \n)$, which we shall refer to as a ``conditioning" below. 
Let us now define
\beaa
\label{eta}
\eta^y_t:=\dfrac{p(t,\xi_t;T,y)}{p(0,\xi_0;T,y)}, \qq t<T, \hQ^0\as,
\eeaa
where $p(\cd, \cd;\cd,\cd)$ is the transition density of $\xi$ under $\hQ^0$, and the process
\beaa
\label{L}
	L^\n_t := \int_{\hR^2}\eta^y_t \nu(dy), \qq t\in [0,T).
\eeaa
 
\begin{defn}
\label{adcond}
The conditioning triplet $(T, \xi_T, \n)$ is called ``proper" if 

(i) supp$(\n) \subseteq$ supp$(\hQ^0\circ \xi_T^{-1})$; and 

(ii) there exist   constants $C, \l>0$, such that 
\bea
\label{etaest}
 0<\sup_{t\in[0, T)}(T-t)\eta^y_t \le{CT}e^{ \frac{\lambda|\xi_0 - y|^2}{T}}, y \in  \hR^2; \q \mbox{\rm and} \q \int_{\hR^2} e
 ^{ \frac{\lambda|\xi_0 - y|^2}{T}}\n(dy)<\infty.
  \eea
\end{defn}

We note that the condition (i) above is relatively easier to verify. In particular, it would be trivial when the diffusion $\xi$ has positive density at time $T$. For condition (ii), we note that $p(s,y;t,x)$ is the fundamental solution to the Kolmogorov backward (parabolic) PDE, then it is well-known that (see, e.g., \cite{Aronson-67, Aronson-68}), for  some constant $c_1$, $c_2$, $\lambda$, $\Lambda >0$, it holds that 
\begin{equation*}
	0<\dfrac{c_1}{t-s}e^{ -\frac{\lambda|y-x|^2}{t-s}}\leq p(s,y;t,x) \leq \dfrac{c_2}{t-s}e^{
	-\frac{\Lambda|y-x|^2}{4(t-s)}}, \q 0\le s<t<T, \q x, y\in\hR^2,
\end{equation*}
Consequently 
 we see that, 
\bea
\label{etabdd}
0<\eta_t \leq\dfrac{c_2T}{c_1(T-t)} e^{-\frac{\Lambda|\xi_t-y|^2}{4(T-t)} + \frac{\lambda|\xi_0 - y|^2}{T}}\leq \dfrac{c_2T}{c_1(T-t)} e^{\frac{\lambda|\xi_0 - y|^2}{T}}, 
\qq t\in[0,T), 
\eea
which leads to the first inequality in (\ref{etaest}). Thus the requirement for the conditioning being ``proper" means that $L^\n_t<\infty$ for all $t\in [0,T)$, $\hQ^0$-a.s..

The following proposition contains some results similar to those  in \cite{FB02}, extended to the 2-dimensional case   but with slightly different assumptions (see also, \cite{FPY, FI93}). Although some proofs are essentially the same, we give a detailed sketch for completeness.

 \begin{prop}
	\label{P1}
	Assume Assumption \ref{assump1}.  Let $(T, \xi_T, \n)$ be a given conditioning. Then, 
		
(i) there exists a unique $\hP^\n\in\sP(\O)$, such that   $\hP^\n\circ \xi_T^{-1}=\n$, and for any $t<T$, any bounded  $X\in \hL^0(\cF_t;\hR^2)$,  it holds that
\bea
\label{EXZT}
\hE^{\hQ^0}\big[X|\xi_T= y\big]=\hE^{\hQ^0}\big[\eta^y_t X\big], \qq t < T, ~y \in \hR^2,\, \hQ^0_{\xi_T}\mbox{-a.s.};
\eea
 
 \ms
(ii) assuming further that $(T, \xi_T, \n)$ is proper, then for any $t<T$,  it holds that

\bea
\label{dPmdQ}
\dfrac{d\hP^{\nu}}{d\hQ^0}\bigg| _{\cF_t}=\int_{\hR^d}\eta_t^y\nu(dy);
\eea

(iii) $L^\n$  is a non-negative  $\hQ^0$-martingale on $[0,T)$, and $\dis L^\n_T:=\lim_{t\to T}L^\n_t$ exists, with $\hE^{\hQ^0}[L^\n_T]\le1$.
\end{prop}

{\it Proof.} 
For the given conditioning $(T, \xi_T, \n)$, let $\hQ^y(\cd)\in \sP(\O)$ be the regular conditional probability defined by 
$\hQ^y(A):=\hQ^0(A|\xi_T=y)$, $A\in\cF_T$, $y\in \hR^2$, and define 
\bea
\label{Pnu}
\dis \hP^\n(A):=\int_{\hR^2} \hQ^y(A)\n(dy), \qq A\in \cF_T.
\eea

We now check (i). That $\hP^\n\circ \xi^{-1}_T=\n$ is obvious. To see  (\ref{EXZT}), 
we define a finite measure on $(\hR^2, \sB(\hR^2))$ by
$\mu^{X |\xi_T}(A):=\int_{\xi_T\in A}X(\o) \hQ^0(d\o)$, $A\in\sB(\hR^2)$. Then, by definition we can write   
\bea
\label{etay}
\mu^{X |\xi_T}(A)=\int_A \hE^{\hQ^0}[X |\xi_T=y] \hQ^0_{\xi_T}(dy)= \int_A \hE^{\hQ^0}[X|\xi_T=y] p(0,z_0; T, y)dy, \q A\in\sB(\hR^2).
\eea
But   since $X\in \hL^0(\cF_t;\hR^2)$, using the Markov property on $\xi$ and Fubini theorem we also have
\bea
\label{etay1}
\mu^{X |\xi_T}(A)
	&=&\int_\O\hE^{\hQ^0}[{\bf 1}_{\{\xi_T\in A\}}X |\cF_t](\o) \hQ^0(d\o)=\int_\O\neg\Big[\neg\int_A \neg p(t, \xi_t(\o); T, y)dy \Big]
	X(\o)\hQ^0(d\o)\nonumber\\
	&=&\int_A\hE^{\hQ^0}[p(t, \xi_t; T, y)X]dy, \q A\in\sB(\hR^2).
\eea
Comparing (\ref{etay}) and (\ref{etay1}), we deduce  (\ref{EXZT}).

(ii) To see (\ref{dPmdQ}), it suffices to show that for any $Z\in \hL^1_{\cF_t}(\hR^2, \hQ^0)$, $t\in[0,T)$, it holds that
\bea
\label{EZ}
\hE^{\hP^\n}[Z]=\hE^{\hQ^0}[L^\n_t Z]=\hE^{\hQ^0}\Big[\int_{\hR^2}\eta^y_t\n(dy)Z\Big].
\eea
By a standard truncation argument, we may assume that $Z$ is bounded. Then by (\ref{EXZT}), we have $\hE^{\hQ^0}[\eta^y_t Z]=\hE^{\hQ^0}[Z|\xi_T=y]=\hE^{\hP^\nu}[Z | \xi_T = y]$, thanks to
definition (\ref{Pnu}), we can write 
\beaa
\hE^{\hP^\n}[Z]=\int_{\hR^2}\hE^{\hQ^0}[Z|\xi_T=y]\n(dy)=\int_{\hR^2}\hE^{\hQ^0}[\eta^y_tZ]\n(dy).
\eeaa
Comparing above to (\ref{EZ}), we see that it suffices to show that $\int_{\hR^2}\hE^{\hQ^0}[|\eta^y_tZ|]\n(dy)<\infty$, so that the Fubini theorem can be applied. But this clearly follows from the boundedness of $Z$ and the assumption that the conditioning   is proper. 

\ms
(iii) Finally, by (\ref{dPmdQ}), $\frac{d\hP^{\nu}}{d\hQ}\Big|_{\scF_t}:=L^\n_t$, $t<T$. Thus,  $L^\n$ is a non-negative $\hQ$-martingale on $[0,T)$.
Furthermore, since $L^\n_t>0$, $t\in[0,T)$, by martingale convergence theorem, $L_T:=\lim_{t\to T}L_t$ exists, and by Fatou's lemma, one easily shows that
$\hE[L^\n_T]\le \lim_{t\to T}\hE[L^\n_t]=1$.
\qed

\begin{rem}
\label{remark3.4}
{\rm (1) The probability  $\hP^\n$ in Proposition \ref{P1} is called the {\it minimal probability}  given the proper conditioning $(T, \xi_T, \n)$. Moreover, Proposition \ref{P1} shows that the assumption (A1) in \cite{FB02} is automatically satisfied in our setting.

(2) Proposition \ref{P1}-(ii) indicates that $\hP^\n$ is absolutely continuous with respect to $\hQ^0$ on each $\cF_t$, 
$0\le t<T$, with the Radon-Nikod\'ym derivative defined by (\ref{dPmdQ}). 
But it does not imply that $\hP^\n$ and $\hQ$ are equivalent 
on $\cF_{t}$, for $t<T$,  neither does it imply that $\hP^\n<\neg\neg<\hQ$ on $ \cF_T$. 
\qed}
\end{rem}

We now turn our attention to a specific conditioning $(T, \xi_T, \n)$ that will lead to the solution to an STPBVP (\ref{vx}).
For notational convenience we shall now simply denote $\xi=(V, X)$, when there is no danger of confusion.  Let 
 $m^*\in \sP(\hR)$ be a law of the underlying asset $V_T$ that is known to the insider. For technical reasons we shall assume   that $m^*$ satisfies the following condition:
 \begin{assum}
 \label{assump2}
 There exists $\l_0>0$ sufficiently large, such that 
 \bea
\label{m*}
\int_{\hR} e^{\l_0 v^2}m^*(dv)<\infty.
\eea
 \end{assum}
 We remark that the Assumption \ref{assump2} is actually not over restrictive. In fact, in light of the well-known Fernique Theorem (cf. \cite{Fer}) (\ref{m*}) covers a large class of normal random variables. Now let us define a probability measure $\nu\in \sP(\hR^2)$ by
\bea
\label{nu}
\nu(A) = \int_{\hR}\mathbbm{1}_A(v,g^{-1}(v))m^*(dv)=\int_{(v,g^{-1}(v))\in A}m^*(dv). 
\eea
That is, the measure $\nu$ concentrates on the graph of the function $v=g^{-1}(x)$, or equivalently $x=g(v)$, thanks to Assumption \ref{assump1}-(iii). 
Furthermore, 
we have the following lemma.
 \begin{lem}
\label{lemma0}
Assume Assumptions \ref{assump1}, \ref{assump2} are in force, with $\l_0$ in (\ref{m*}) being sufficiently large. 
Let $\xi$ be the solution to (\ref{SDEVXQ}), and $\n\in \sP(\hR^2)$ be defined by (\ref{nu}). Then, $(T, \xi_T, \n)$ is a proper conditioning. Furthermore, if $\hP^\n$ is the minimum probability given $(T, \xi_T, \n)$, then it holds that
\bea
\label{PnuT}
\hP^\n\{V_T= g(X_T)\}=1.
\eea
\end{lem}

{\it Proof.} Since under Assumption \ref{assump1} $\xi$ is a diffusion process with  positive  transition density function (cf. e.g., \cite{FW}), we have supp$(\hQ\circ \xi_T^{-1})=\hR^2$. Furthermore, by definition of $\n$ (\ref{nu}), for the constants $ \l>0$ in (\ref{etabdd}) we deduce from (\ref{m*}) that
\beaa
\int_{\hR^2} e^{\frac{\lambda|\xi_0 - y|^2}{T}}\n(dy)=\int_{\hR} e^{\frac{\lambda[(v_0 - v)^2+(g(x_0)-v)^2]}{T}}m^*(dv)<\infty,
\eeaa
provided that $\l_0\ge \frac{2\l}{T}$, where $\l_0$ is the constant in Assumption \ref{assump2}. That is, $(T, \xi_T, \n)$ is proper. 

To show the second assertion, first note that $g$ is strictly increasing, the graphs of $g$ and $g^{-1}$, as the subset of $\hR^2$, are identical. Let us denote
$ \G:=\{(g(x),x): x\in \hR\}=\{(v, g^{-1}(v)):v\in\hR\}\subseteq \hR^2$.
Then, by definition (\ref{nu}) we see that $\n(A)=1$ if and only if $\G\subseteq A$. In particular, $\n(\G)=1$. Consequently, by definition of the minimum probability, we have 
$$\hP^\n(V_T=g(X_T)\}=\hP^\n \circ \xi_T^{-1}(\G) =\nu(\G)=1,
$$ 
proving (\ref{PnuT}).
\qed

\begin{rem}
\label{remark3.6} 
{\rm (1) Given (\ref{PnuT}), and the fact that  $\xi=(V, X)$ has continuous paths under $\hQ^0$, it is readily seen that 
$\hP^\nu \{\lim_{t\to T} V _t = V_T=g(X_T) = \lim_{t\to T} g(X_t)\}=1$. 
Consequently, we have 
\bea
\label{VgXT}
	\hP^\nu \{\lim_{t\to T} V_t = \lim_{t\to T} g(X_t)\}=1. 
\eea
This, together with Proposition \ref{prop37}, indicates that as far as the solution to the two-point boundary value problem is concerned, without the specific requirement of Markovian bridge, the SDE (\ref{SDEVXPn}) would be a desirable candidate, except for a slight difference on the drift coefficients.

\ss
(2) By Proposition \ref{P1}-(iii), $L^\n$ is a closeable supermartingale on $[0, T]$. But  it cannot be a martingale, unless $\hQ^0\{V_T=g(X_T)\}=1$, which is obviously not true in general.  Thus $\hP^\n$ cannot be absolutely continuous with respect to $\hQ^0$ on $\cF_T$, 
as we pointed out in Remark \ref{remark3.4}.  
\qed
}
\end{rem}

To end this section,  let us define, for any proper conditioning $(T, \xi_T, \n)$, a function
\bea
\label{phi}
\varphi(t, z)=\int_{\hR^2}\dfrac{p(t,z;T,y)}{p(0,z_0;T,y)}\nu(dy), \qq z=(v,x), 
\eea
where $p$ is the transition density of  $\xi$ under $\hQ^0$. Then, clearly, $\varphi(0,z_0)=1$, and  $L_t=L^\n_t=\varphi(t, \xi_t)$, $t\in [0,T)$. Now, applying It\^o's formula we have 
\bea
\label{Itovarf}
 L_t=\varphi(t,\xi_t) =1 + \int_0^t [\pa_t \varphi(s,\xi_s)+ \scL[\varphi](s,\xi_s)]ds + \int_0^t \big(\nabla \varphi(s,\xi_s), \bar{\si}(s, \xi_s)dB^0_s\big), 
\eea
where 
$\scL[\varphi](t,z):= (\bar{b}, \nabla\varphi)(t,z)+ \tr[D^2\varphi\bar{\si}\bar{\si}^T](t,z)$, and $\bar b :=(b,  \m)^T$, $\bar{\si}:=\mbox{\rm diag}[\si,\rho]$. Since  by Proposition \ref{P1}-(iii),
$L$ is a $\hQ^0$-martingale for $t\in[0, T)$, we conclude that $\varphi$ must satisfy the following PDE
(noting the definition of $\bar{b}$ and $\bar{\si}$) 
for $t\in[0,T)$ and $z=(v,x)\in\hR^2$, 
\begin{equation}
\label{pdephi}
\left\{\ba{lll}
\dfrac{\pa \varphi(t,z)}{\pa t}+b(t,z)\dfrac{\pa \varphi}{\pa v}+\mu(t,x)\dfrac{\pa \varphi(t,z)}{\pa x}+\dfrac{1}{2}\sigma^2(t,z)\dfrac{\pa ^2 \varphi}{\pa v^2}+\dfrac{1}{2}\rho^2(t,x)\dfrac{\pa^2 \varphi(t,z)}{\pa x^2}=0;\ms\\
\varphi(0, v_0,x_0)=1.
\ea\right.
\end{equation}
Consequently, it follows from (\ref{Itovarf}) that  
\bea
\label{L0}
	dL_t=d\varphi(t, \xi_t)=  \big(\nabla \varphi(t,\xi_t), \bar{\si}(t, \xi_t)d B^0_t\big)=L_t(\th_t,  dB^0_t), \qq L_0=1, \q t\in [0, T),	
\eea
where  $\th_t :=\bar{\si}^T(t, \xi_t)\frac{\nabla\varphi(t,\xi_t)}{\varphi(t,\xi_t)}=\bar{\si}^T(t, \xi_t)\nabla [\ln \varphi(t, \xi_t)]$, $t\in [0,T)$.
 Denote 
$ W_t = B^0_t - \int_0^t \th_s ds$, then by Girsanov theorem, $\{W_t\}$ is a 2-dimensional $\hP^\n$-Brownian motion on $[0,T)$. 
In other words, we have proved the following 2-dimensional extension of a result in \cite{FB02}.
\begin{prop}[{\cite[Proposition 37]{FB02}}]
\label{prop37} 
Assume Assumption \ref{assump1}, and let  $\hP^\n$ be the minimal probability corresponding to the conditioning $(T, \xi_T, \n)$, where $\xi = (V,X)$ is the strong solution to \eqref{SDEVXQ}. Then, under $\hP^\n$, $\xi$ solves  the following SDE:
\bea
\label{SDEVXPn}
d\xi_t=[\bar{b}(t, \xi_t)+\bar{\si}(t, \xi_t)\th_t]dt +\bar{\si}(t, \xi_t)dW_t= \hat{b}(t, \xi_t)dt +\bar{\si}(t, \xi_t) dW_t, \q \xi_0  =z,
\eea
where $(\bar{b}, \bar{\si})$ are the same as those in  (\ref{Itovarf}),  $\hat{b}(t, v,x):=\bar{b}(t, v,x)+\bar{\si}\bar{\si}^T(t,v,x)\nabla[\ln\varphi(t,v,x)]$, and  $\th_t :=(\theta^1_t, \theta^2_t)^T=\bar{\si}^T(t, \xi_t)\nabla[\ln\varphi(t, \xi_t)]=\frac1\varphi ( \pa_v \varphi \sigma, \pa_x \varphi \rho)^T(t, \xi_t)$; 
 $\varphi$ is defined by (\ref{phi}); and $W=(W^1, W^2)$ is a $\hP^\n$-Brownian motion. 
\end{prop}

 \section{A Stochastic Two-Point Boundary Value Problem}
 \setcounter{equation}{0}
 
We are now ready to study the STPBVP (\ref{vx}) and compare it to the  well-known {\it dynamic Markov bridge} in the literature.  We begin by giving the precise definition of the STPBVP. 
\begin{defn}
\label{TPBVP}
A six-tuple   $(\hP, B^1, B^2, V,X, \a)$ is called a (weak) solution of a stochastic Two-Point Boundary Value Problem (STPBVP) 
on $[0, T]$ if 
\ss

(i) $\hP\in \sP(\O)$ and $B=(B^1, B^2)$ is a $\hP$-Brownian motion on $[0,T]$;

\ss
(ii) $\a\in \sU_{ad}$, and $(V,X, \a)$ satisfies the SDE  on $(\O, \cF, \hP)$: 
\bea
\label{SDEVX}
	\left\{\begin{array}{lll}
	dV_t = b(t, V_t, X_t)dt + \sigma(t, V_t, X_t) dB^1_t, \qq & V_0=v; \ms\\
	dX_t = \big(\mu(t, X_t)+\a_t\rho(t,X_t)\big)dt +\rho(t,X_t)dB^2_t, & X_0=x,
	\end{array}\right.
	 \qq t\in [0,T), \hP\as;
\eea

(iii) $\lim_{t\nearrow T}[V_t-g(X_t)]=0$, $\hP$-a.s.;

\ss
In particular, $(V,X, \a)$ is called the solution to a Markovian STPBVP, if $\a_t=u(t, V_t, X_t)$, $t\in [0,T)$, for some measurable function $u$, and $(V,X)$ is an $\hF^{V,X}$-Markov process on $[0,T)$.
\qed
\end{defn}

\begin{rem}
\label{remark3a}
{\rm 
(i) For notational clarity we shall often refer to (\ref{SDEVX}) as a ``STPBVP$(b, \si, \m, \rho)$" when necessary, and write the solution $(V, X, \a)$ to a STPBVP as $(V^\a, X^\a)$ for convenience. 

\ss
(ii) Comparing Definition \ref{TPBVP} to that of a dynamic Markov bridge (see, e.g., \cite{CCD11}), we see that, if the coefficients
$b$ and $\si$ are independent of $X$ and $\m\equiv 0$, then a Markovian TPBVP is essentially a dynamic Markov bridge without requiring that $X$ be a local martingale with respect to its own filtration $\hF^X$. Such a difference will make the results of this paper and those in the existing literature mutually exclusive. 
\qed}
\end{rem}

To construct a weak solution, we first recall (\ref{L0}) and the $\hP^\nu$-Brownian motion  $W_t= B^0_t - \int_0^t \theta_s ds$; $t\in[0, T)$,
where $\th_t :=(\theta^1_t, \theta^2_t)^T=\bar{\si}^T(t, \xi_t)\nabla[\ln\varphi(t, \xi_t)]=\frac1\varphi ( \pa_v \varphi \sigma, \pa_x \varphi \rho)^T(t, \xi_t)$, $t\in[0,T)$,  and under  $\hP^\n$ the process $\xi_t:=(V_t, X_t)^T$ satisfies the SDE (\ref{SDEVXPn}).
We note that although the coefficient $\hat{b}$ in (\ref{SDEVXPn}) is explicitly defined, it depends on the solution of an ill-posed parabolic PDE (\ref{pdephi}), its behavior is a bit hard to analyze. The following lemma is useful to note. 
\begin{lem}
\label{lemma1}
Let $(T, \xi_T, \n)$ be the conditioning in Lemma \ref{lemma0}, and $\hP^\n$  the corresponding minimum probability. Then, 
it holds that $\hL^p_{\cF_t}(\hR^d;\hQ^0)\subset \hL^p_{\cF_t}(\hR^d; \hP^\n)$, $t<T$. Specifically, for any $T_0<T$, there exists a constant $C_{T_0}>0$, that depends only on the coefficients $(b, \si, \m, \rho)$, and $T_0$,  such that, for any $X\in \cF_t$, $t\in[0,T_0]$, it holds that
\bea
\label{PnbQ}
 \hE^{\hP^\n}[|X|^p]\le C_{T_0} \hE^{\hQ^0}[|X|^p].
 \eea
 In particular, the $\hQ^0$-diffusion process $\xi$
is well-defined for $t\in [0,T)$ on the probability space $(\O, \cF, \hP^\n)$, and $\hP^\n\{\int_0^{T_0}|\xi_t|^2dt<\infty\}=1$, for any $T_0<T$. 
\end{lem}

{\it Proof.} We first note that given $T_0<T$, and $X\in \cF_t$, $t\le T_0$, by Lemma \ref{lemma0}-(ii)
\beaa
\hE^{\hP^\n}[|X|^p] = \hE^{\hQ^0}[L_{T_0}|X|^p]\le C_{T_0}  \hE^{\hQ^0}[|X|^p],
\eeaa
where  $C_{T_0}:=\frac{\widetilde{C}T}{T-T_0}\int_{\hR^2}e^{\frac{\l|\xi_0-y|^2}{T}}\n(dy)$, proving (\ref{PnbQ}). 
The rest of the proof is obvious.
\qed

Now for any $n\in\hN$, define the stopping time 
$\tau_n := \inf \{ t>0 : |\theta_t |\geq n\} \wedge T$, and then denote
$\theta^{(n)}_t:=\th_{t\wedge \t_n}$, $t\in[0,T]$.
Clearly,  under probability $\hP^\nu$, for each $n\in\hN$, the SDE
 \bea
\label{Zn}
d\xi^{(n)}_t=[\bar{b}(t, \xi^{(n)}_t)+\bar{\si}(t, \xi^{(n)}_t)\th^{(n)}_t]dt +\bar{\si}(t, \xi^{(n)}_t)dW_t,
 \q \xi^{(n)}_0  =z,
\eea
 is well-posed on $[0,T]$, and $\xi^{(n)}_t\equiv \xi_t$, $t<\tau_n$. Furthermore, 
since $\tau_n \nearrow T$ as $n \to \infty$,   we can simply define $\xi_t = \xi_t^{(n)}$, $t\in[0,\tau_n]$, $n=1,2,\cds$. Then it holds that $\lim_{t\nearrow T} \xi^{(n)}_t = 
\lim_{t\nearrow T}\xi_t = \xi_T$. 
In particular, we have
$\hP^\nu \{ \lim_{t\nearrow T} V_t = \lim_{t\nearrow T} g(X_t)\} =1$.

We now write $\th^{(n)}_t=(\th^{1,n}_t, \th^{2,n }_t)$, $t\in[0,T]$. 
Since    $\theta^{1,n}_t$ is bounded by $n$, and $\theta^{1,n}_t = \theta^{1,n+1}_t$, on $[0,\t_n]$.
By Girsanov's theorem, there exists a  family of probabilities $\{\bar{\hP}^{(n)}\}_{n\geq 1}$ on $(\O, \cF)$ by
\begin{equation*}
	\dfrac{d\bar{\hP}^{(n)}}{d\hP^\nu} \Big| _{\scF_T}=\sE( \theta^{1,n}_T):=\exp\Big\{\int_0^T\th^{1,n}_sdW^1_s-\frac12\int_0^T|\th^{1,n}_s|^2ds\Big\}.
\end{equation*}
Then for each $n\in\hN$, the process  $\bar{B}^{(n)}_t = (\bar{B}^{1,n}_t, W^2_t) := (W^1_t - \int_0^t\theta^{1,n}_s ds, W^2_t)$, $t\in [0,T]$,
is a 2-dim $\bar{\hP}^{(n)}$-Brownian motion. Moreover, by the property of $\{\th^n\}$, we must have 
\bea
\label{bPn}
\dfrac{d\bar{\hP}^{(n+1)}}{d\hP^\nu} \Big| _{\cF_{\t_n}} = \sE(\theta^{1,n+1}_{\t_n})  = \sE(\theta^{1,n}_{\t_n}) = \dfrac{d\bar{\hP}^{(n)}}{d\hP^\nu} \Big| _{\cF_{\t_n}}. 
\eea
Consequently, we have 
$	\bar{\hP}^{(n+1)} \big| _{\cF_{\t_n}} = \bar{\hP}^{(n)}\big| _{\cF_{\t_n}}$, 	
and $\bar{B}^{(n+1)}_t =\bar{B}^{(n)}_t$, $ t \in[0,\t_n]$, for each $n\in\hN$. 
Observing that $\t_n \nearrow T$ as $n \to \infty$, we can define a new probability measure $\bar\hP$  on $(\O, \cF_{T-})$  by 
\bea
\label{bP}
\bar{\hP}|_{\cF_{\t_n}}:=\bar{\hP}^{(n)}|_{\cF_{\t_n}}, \qq n\in\hN,
\eea
then $\bar{\hP}<\neg\neg< \hP^\n$ on $\cF_t$, $t\in [0, T)$.  Furthermore, if we define  $\bar{B}_t = \bar{B}^{(n)}_t$, $t\in [0,\t_n]$, $n\in\hN$, then $\bar B$ is a $\bar\hP$-Brownian motion on $[0,T)$, and under $\bar\hP$, the process $\xi=(V,X)$ satisfies the  SDE:
\bea
\label{SDEVXpbar}
\left\{\begin{array}{lll}
	dV_t = b(t,V_t ,X_t )dt+\sigma(t,V_t ,X_t) d\bar{B}^1_t, \qq\qq\qq & V_0=v; \ms\\
	dX_t= \big(\mu(t,X_t )+\rho(t,X_t) \theta^2_t\big)dt +\rho(t,X_t )dW^2_t ,& X_0 =x;
\end{array}\right. \q t\in [0,T).
\eea
Comparing (\ref{SDEVXpbar}) and (\ref{SDEVX}) and noting the facts (\ref{VgXT}) and $\bar{\hP}|_{\cF_t}<\neg\neg< \hP^\n|_{\cF_t}$, $t\in [0, T)$, we see that $( \bar\hP,  \bar{B},  V, X, \th^2)$ should be a weak solution to STPBVP (\ref{SDEVX}). We have the following result.
\begin{prop}
\label{limitbarP}
	Assume Assumption \ref{assump1}. Then there exists a weak solution to  STPBVP (\ref{SDEVX}). 
Furthermore,  if  $( \hP, B, V , X , \a)$ denotes the weak solution, then 
$\hP$ can be chosen so that $\hP_{\cF_t}<\neg\neg<\hQ^0|_{\cF_t}$, $t<T$, and denoting $V_T:=V_{T-}=\lim_{t\nearrow T}V_t$, it holds that $\hP\circ (V_{T})^{-1} =m^*$. 
\end{prop}

{\it Proof.}  Consider the probability $\bar\hP$ defined by (\ref{bPn}), (\ref{bP}) and SDE (\ref{SDEVXpbar}). We first claim  $\bar{\hP}<\neg\neg<\hP^\nu$ on $\cF_{T-}$. Indeed,  let
$\sA:=\{\cG\subset\cF:  \bar{\hP}<\neg\neg < \hP^\nu \mbox{~on $\cG$}\}$, then $\cF_{\t_n}\in\sA$, $n\in \hN$. Since $\t_n\nearrow T$, we  have $\cF_{T-}=\bigvee_n\cF_{\t_n}$ (see, e.g., \cite{PR}), and thus
 $\cF_{T-} \in \sA$, thanks to the  Monotone Class Theorem.
 	
Next, since
$\{ \lim_{t\nearrow T} V_t \neq \lim_{t\nearrow T} g(X_t)\} = \bigcup_{m }\bigcap_{N}\bigcup_{r\in {\bf Q}(T-\frac1N, T)}\big\{ |V_{r} - g(X_{r})|\ge\frac1m\big\}\in \cF_{T-}$, where $\bf Q$ is the rationals in $\hR_+$, and ${\bf Q}(A):={\bf Q}\cap A$, $A\in \cB(\hR)$, 
and $\bar{\hP}<\neg\neg < \hP^\nu$ on $\cF_{T-}$, we have
$\bar{\hP} \{ \lim_{t\nearrow T} V_t \neq \lim_{t\nearrow T} g(X_t)\} = 0$, thanks to (\ref{VgXT}). That is, $\bar{\hP} 
\{ \lim_{t\nearrow T} V_t = \lim_{t\nearrow T} g(X_t)\} =1$.
Now let 
$\a=\th^2$ in SDE (\ref{SDEVXpbar}), we see that  $( \bar\hP,  \bar{B},  V, X, \a)$ is a weak solution to STPBVP (\ref{SDEVX}). 

It remains to check the last statement.  To this end, let $\xi=(V,X)$. Since $\bar\hP<\neg\neg<\hP^\n<\neg\neg<\hQ^0$ on $\cF_{T-}$ and 
$\hQ^0\{\xi\in \hC([0, T];\hR^2)\} =1$, we can naturally extend  $\xi$ to $[0,T]$ by  setting $\xi_T=\lim_{t\nearrow T}\xi_t$ so that ${\hP}^\n\{\xi\in \hC([0, T];\hR^2)\} =\bar{\hP}\{\xi\in \hC([0, T];\hR^2)\}=1$ as well. We first claim that $\hP^\nu\circ V_T^{-1}=m^*$. Indeed, let $B \in\sB( \hR)$ and $A:=B\times \hR\in \sB(\hR^2)$. By (\ref{nu}) we have
$B = \{v : (v,g^{-1}(v))\in A\}$, and 
$\hP^\nu\{V_T \in B\}= \hP^\nu\{ (V_T,X_T)\in A\} = \nu\{A\} = m^*\{B\}$. That is, $\hP^\nu\circ V_T^{-1}=m^*$.

To see $\bar\hP\circ V^{-1}_T=m^*$, we note that $\xi=(V, X)$ is the unique strong solution to SDE (\ref{SDEVXQ}) under 
$\hQ^0$ with canonical process $B^0=(B^1, Y)$. Therefore we can write $\xi_t(\o)=\Phi(t, B^0_{\cd\wedge t}(\o))=\Phi(t, \o)$, $(t, \o)\in[0,T]\times\O$, for some (progressively) measurable function $\Phi: [0,T]\times\O\mapsto \hR^2$. Consequently, we can write 
$\th^2_t(\o) =\pa_x \ln\varphi(t, \xi_t(\o))=\pa_x \ln\varphi(t, \Phi(t, \o))$, $(t,\o)\in[0,T]\times \O$.
By virtue of Lemma \ref{lemma1}, the process $\th^2$ is  well-defined on $[0,T)\times\O$, $\bar\hP$-a.s. and $\th^2_t\in \hL^2(\bar\hP)$, for $t\in[0,T)$. 

Now let us denote the solutions to \eqref{SDEVXPn} and \eqref{SDEVXpbar}  as $(\tilde{V}_t,\tilde{X}_t)$ and $(\bar{V}_t, \bar{X}_t)$ respectively. Then we see that $((\tilde{X}_t, W^2_t), \hP^\nu)$ and $((\bar{X}_t, W^2_t), \bar{\hP})$ are two weak solutions to the same SDE, well-defined on any $[0,T_0]\subset[0,T)$. 
Consequently, we have $\hP^\n\circ {\tilde X}^{-1}=\bar{\hP}\circ \bar{X}^{-1}$ on $[0, T_0]$ for any $T_0<T$. Extending the solution to $[0,T]$, we have 
${\hP^\n} \circ \tilde{X}_T ^{-1} = \bar{\hP} \circ \bar{X}_T ^{-1} $. Since 
$V_T = g(X_T)$, both $\bar{\hP}$-a.s. and $\hP^\n$-a.s., we obtain that $\bar{\hP} \circ V_T ^{-1} = \hP^\n\circ V_T^{-1}=m^*$,  proving the proposition. 
 \qed

\ss

{\bf Uniqueness in law.} Let us now turn to the issue of uniqueness. 
To begin with let us recall 
that the weak solution $(\bar{\hP}, \bar{B}, V, X, \a)$ that 
we constructed has the following properties:

\ss
(i) there exists a sequence of $\bar{\hP}$-stopping times $\{\t_n\}$, and a sequence of probabilities $\bar{\hP}^{(n)}$ on $(\O, \cF)$, such 
that $\t_n\nearrow T$, $\bar\hP$-a.s., and 
$\bar{\hP}|_{\cF_{\t_n}}=\bar\hP^{(n)}|_{\cF_{\t_n}}$, $n\in\hN$;

(ii) for each $n\in\hN$,  $\bar{B}=\bar{B}^{(n)}$ on $[0,\t_n]$, where $\bar{B}^{(n)}= (\bar{B}^{(n,1)}, \bar{B}^{(n,2)})$ is a $\hP^{(n)}$-Brownian motion on
$[0,T]$; 

(iii)  the solution $(\bar{V}, \bar{X})=(V^{(n)}, X^{(n)})$ on $[0,\t_n]$, where $(V^{(n)}, X^{(n)})$ is a (pathwisely) unique solution to the following SDE, defined on $[0,T]$:
\bea
\label{SDEVXpbar1}
\left\{\begin{array}{lll}
	dV_t = b(t,V_t ,X_t )dt+\sigma(t,V_t ,X_t) dB^{(n,1)}_t, \qq & V_0=v; \ms\\
	dX_t= \big(\mu(t,X_t )+\rho(t,X_t) \a^{(n)}_t\big)dt +\rho(t,X_t )dB^{(n,2)}_t, & X_0 =x,
\end{array}\right.
\eea
where $|\a^{(n)}_t|\le M_n$, $t\in[0,T]$, for some $M_n>0$; and $\a^{(n+1)}_t=\a^{(n)}_t$,  $t\in[0,\t_n]$, $\bar\hP$-a.s.;

\ms
(iv) $\bar{\hP}|_{\cF_t}<\neg\neg<\hP^{\n}|_{\cF_t}<\neg\neg<\hQ^0|_{\cF_t}$, $t\in [0,T)$. 

\ms
In what follows we shall 
denote $(\bar{\hP}, \{\t_n\})$ to specify that $\bar{\hP}$ is ``announced" by $\{\t_n\}$, and make use of the following definitions in the spirit of the so-called ``{\it $\hQ^0$-weak solutions}" in \cite{MSZ}.
\begin{defn}
\label{nested}
We call a weak solution $(\bar{\hP}, \bar{V}, \bar{X}, \bar{B},  \a)$ of STPBVP (\ref{TPBVP}) satisfying (i)--(iii) above a ``{\it nested weak solution}" and the corresponding family of stopping times $\{\t_n\}$ the ``{\it announcing sequence}" of probability $\bar\hP$. 
We  call $(\{\t_n\}, \a)$ the characteristic pair of the  weak solution. 

Furthermore, a nested weak solution is called a $\hP^\n$-weak solution if it also satisfies (iv).
\qed
\end{defn}

\begin{rem}
\label{remark3.3}
{\rm
Comparing to the usual SDEs, the characteristic pair $(\{\t_n\}, \a)$ is important in determining a  solution to an STPBVP. Note that if $\{\t_n^1\}, \{\t_n^2\}$ are two announcing sequences of stopping times, then so is
$\{\t_n^1\wedge \t_n^2\}$. Thus the  weak solution is independent of  the choice of the announcing sequence $\{\t_n\}$.
Since  the process $\a$ determines the coefficient of SDE (\ref{SDEVXpbar}), whence the solution, we often specify its role by calling $(\bar{\hP}, \bar{V}, \bar{X}, \bar{B},  \a)$ the $\a$-weak solution.
 \qed}
\end{rem}

\begin{defn}
\label{Def1}
We say that the pathwise uniqueness holds for  STPBVP (\ref{TPBVP}), if for two nested solutions $({\hP}^i, \xi^i=( {V}^i,  {X}^i),  {B}^i,  {\a}^i)$, $i=1,2$ of  (\ref{TPBVP}) on $[0,T)$, such that $\hP^1=\hP^2=\hP$, $\xi^1_0=\xi^2_0$, and $\hP\{\a^1_t=\a^2_t, 
~ B^1_t=B^2_t, ~ t\in[0,T)\}=1$, then ${\hP}\{\xi^1_t= \xi^2_t, ~  t\in [0,T_0]\}=1$, for any $T_0<T$.
\end{defn}
\begin{rem}
\label{pwuniq}
{\rm We note that the time $T_0 $ in Definition \ref{Def1} can be changed to any stopping time $\t$ with $\hP\{\t<T\}=1$. In fact, the following two statements are equivalent: (i) the pathwise uniqueness holds on $[0,T_0]$, for any $T_0 < T$; and
(ii) there exists a sequence of stopping time $\{\tau_n, n\geq 1\}$, $\lim_{n \to \infty} \tau_n = T$ almost surely, such that the pathwise uniqueness holds on $[0,\tau_n]$, for each $n\geq 1$.
Indeed, let $ (\hP^i, \xi^i=(V^i, X^i))$, $i=1, 2$, be two nested solutions as in Definition \ref{Def1}, and denote $\D\xi:=\xi^1_t - \xi^2_t$, then we obtain
	\begin{align*}
		\hE\Big[\sup_{t\in[0,T_0]} |\D\xi|\Big] 
		&\leq
		\hE\Big[\sup_{t\in[0, \tau_n]} \vert \D\xi\vert{\bf{1}}_{\{T_0\leq \tau_n\}} \Big] + \hE\Big[\sup_{t\in[0,T_0]} \vert \D\xi\vert{\bf{1}}_{\{T_0>\tau_n\}}\Big]\le \hE\Big[\sup_{t\in[0,T_0]} \vert \D\xi\vert{\bf{1}}_{\{T_0>\tau_n\}}\Big];
	\end{align*}
and similarly, for any $T_0<T$, 
	\begin{align*}
		\hE\Big[\sup_{t\in[0,\tau]} |\D\xi|\Big] 
		&\leq
		\hE\Big[\sup_{t\in[0, T_0]} \vert \D\xi\vert{\bf{1}}_{\{\tau\leq T_0\}} \Big] + \hE\Big[\sup_{t\in[0,\tau]} \vert \D\xi\vert{\bf{1}}_{\{\tau>T_0\}}\Big]\le \hE\Big[\sup_{t\in[0,\tau]} \vert \D\xi\vert{\bf{1}}_{\{\tau>T_0\}}\Big].
	\end{align*}
Since $\lim_{n\to\infty} \hP \{T_0 > \tau_n\} = 0$ and 	$\lim_{T_0\nearrow T} \hP \{\t> T_0\} = 0$, it is readily seen that the statements (i) and (ii) above  are equivalent, and $T_0$ in Definition \ref{Def1} can be replaced by any $\t <T$.
\qed}
\end{rem}

The definition of the uniqueness in law for the STPBVP is a bit more involved. First note that the component ``$\a$" of the solution 
is part of the drift coefficient of the SDE (\ref{SDEVXpbar}), and in general it is not unique. Thus the uniqueness of the solution, even in the weak sense, depends   on how the process $\a$ is properly fixed. 
To this end,  denote  $\sA:=\{A\in \sB([0,T])\otimes\cF: A_t\in \cF_t, ~t\in[0,T]\}$, where $A_t$ is the $t$-section of $A $; and 
 denote all $\sA$-measurable functions by $\hL^0_\sA([0,T]\times\O)$. 
We should note that the space $\hL^0_{\sA}([0,T]\times\O)$ is independent of any probability measure, and we can therefore use it to identify the $\a$-component of the solution in an ``universal" way.
\begin{defn}
\label{uniq-in-law}
We say that the nested weak solution to the STPBVP (\ref{TPBVP}) is unique in law, if for any two  $\a$-weak solutions 
$(\bar{\hP}^i, \bar{V}^i, \bar{X}^i, \bar{B}^i, \bar{\a}^i )$, $i=1,2$ of  (\ref{TPBVP}) on $[0,T)$, such that  $(v^1,x^1)=(v^2,x^2)$; 
 $\bar{\hP}^1\circ (\t_n^1)^{-1}=\bar{\hP}^2\circ (\t_n^2)^{-1}$, $n\in\hN$; and $\bar{\hP}^i\{\bar{\a}^i_t=\a_t, ~ t\in[0,T)\}=1$,
 $i=1,2$, for some $\a\in \hL^0_{\cA}([0,T]\times\O)$, 
then for any cylindrical set $E^{A_1,\dots, A_n}_{t_1,\cds,t_n}
:=\{ ({\bv,\bx})\in \hC([0,T];\hR^2): (\bv, \bx)(t_i)\in A_i, ~ i=1, \cds, n\}$,  where  $0\le t_1<t_2<\cds <t_n<T$  and $A_i\in\sB(\hR^2)$, $i=1, \cds, n$, it holds that 
\beaa
\label{unique}
\bar{\hP}^1\circ(\bar{V}^1, \bar{X}^1)^{-1}\{E^{A_1,\dots, A_n}_{t_1,\cds,t_n}\}= \bar{\hP}^2\circ(\bar{V}^2, \bar{X}^2)^{-1}\{E^{A_1,\dots, A_n}_{t_1,\cds,t_n}\}.
\eeaa
\end{defn}
We now give the main theorem of this subsection. 
\begin{prop}
\label{uniqinlaw}
Assume  Assumption \ref{assump1}. Then, the Markovian
$\hP^\n$-weak solution to  STPBVP (\ref{TPBVP}) is unique in law. 
\end{prop}
Before we prove Proposition \ref{uniqinlaw}, we first prove a lemma that is interesting in its own right. 
\begin{lem}
\label{alpha}
Assume Assumption \ref{assump1}, and let $(\bar\hP, \bar{\xi}, \bar{\a})$ be a nested Markovian weak solution with $\bar{\a}_t=u(t, \bar{\xi}_t)$, $u\in\hL^0([0,T]\times\hR^2)$, such that  
$\bar{\hP}\{\bar{\a}_t=\a_t, ~t\in[0,T)\}=1$ for some $\a\in \hL^0_{\sA}([0,T]\times \O)$. Then 
$\a_t(\o)=u(t, \Phi(t,\o))$, $dt\otimes d\bar{\hP}$-a.e.-$(t, \o)\in[0,T)\times\O$,  for some $\Phi\in \hL^0_{\sA}([0,T)\times\O)$. \end{lem}

{\it Proof.}
Let $(\bar{\hP}, \bar{\xi}, \bar{\a})$ be the nested Markovian weak solution. Then $\bar{\a}_t =u(t, \bar{\xi}_t)$, $t\in[0,T]$, for
some $u\in\hL^0([0,T]\times\hR^2)$. By Definition \ref{nested}, the solution $\bar{\xi}$ is the pathwisely unique weak solution of SDE (\ref{SDEVXpbar1}) on any $[0,\t_n]$, $n\ge 1$, whence on $[0,T_0]$, for any $T_0<T$, thanks to Remark \ref{pwuniq}. Thus, by Yamada-Watanabe theorem, for any $T_0<T$,  $\bar{\xi}$ is the pathwisely unique strong solution on $[0,T_0]$, and there exists a $\Phi^{T_0}\in
\hL^0_\sA([0,T_0]\times\O)$, such that $\bar{\xi}_t=\Phi^{T_0}(t, \cd)$, $t\in[0,T_0]$, $\bar\hP$-a.s.. As before, we can define a $\Phi\in\hL^0([0,T]\times\O)$ so that $\Phi(t, \cd)=\Phi^{T_n}(t, \cd)$,  $t\in[0, T_n]$, for any sequence $T_n\nearrow T$, and  $\bar{\xi}_t=\Phi(t, \cd)$, $t\in[0,T)$, $\bar\hP$-a.s.. Since
$\bar{\a}_t=u(t, \bar{\xi}_t)=u(t, \Phi(t, \cd))$, $t\in[0,T)$, by assumption, we have $\a_t(\o)=\bar{\a}_t(\o)=u(t, \Phi(t, \o))$,  $dt\otimes d\hP$-a.e., 
proving the lemma.
\qed

\ss
[{\it Proof of Proposition \ref{uniqinlaw}.}] 
 Let $(\bar{\hP}^i, \bar\xi^i_t=(\bar{V}^i, \bar{X}^i), \bar{B}^i, \a^i)$, $i=1,2$, be two Markovian weak solutions of (\ref{TPBVP}) on $[0,T)$, with characteristic pair $(\{\t^i_m\}, \a^i)$, $i=1,2$. Without loss of generality,  we assume that $\{\t^i_m\}$ is the exit time of $\a^i = u(t, \bar{\xi}^i)$, $i =1,2$, from the interval $[-m, m]$.

Next, let the cylindrical set $E^{A_1,\dots, A_n}_{t_1,\dots,t_n}$ be given, with $t_n<T$. Since $\t^i_m\nearrow T$, we can write
\beaa
	(\bar{\xi}^i)^{-1}(E^{A_1,\dots, A_n}_{t_1,\dots,t_n} )
		=\bigcap_{j=1}^n(\bar{\xi}^i_{t_j})^{-1}(A_j)=\bigcup_{m=1}^\infty \bigcap_{j=1}^n \{ \t^i_m  \geq t_j\}\cap (\bar{\xi}^i_{t_j})^{-1}(A_j), \q i=1,2. 
\eeaa
Denoting $E^i_{j,m} := \{ \t^i_m  \geq t_j\}\cap (\bar{\xi}^i_{t_j})^{-1}(A_j)= \{ \t^i_m  \geq t_j\}\cap (\bar{\xi}^{i, (m)}_{t_j})^{-1}(A_j)$, i=1,2, we claim that $E^i_{j,m} \in \cF_{\t^i_m}$, for each $i, j, m$. 
Indeed,  fix $i, j$, and $m$,  for $ t\in[0,T)$, one has
\vspace{-2mm}
	$$\{\t^i_m \leq t\}\cap  E^i_{j,m} = \{ t_j\le \t^i_m \leq t\}\cap (\bar{\xi}^{i, (m)}_{t_j})^{-1}(A_j) \in \cF_t, 
	\q i=1,2.$$
That is, $E^i_{j,m}\in \cF_{\t^i_m}$, whence $\hat E^i_m:=\bigcap_{j=1}^n E^i_{j,m} \in \cF_{\t^i_m}$, $i=1,2$. 

On the other hand, note that the set $\hat E_m$ is increasing in $m$, thanks to the extension nature of solutions $\bar{\xi}^{i,(m)}$. Thus, noting that $\bar{\hP}^i | _{\cF_{\t^i_m} }= \bar{\hP}^{i,(m)}| _{\cF_{\t^i_m}}$, for $i=1,2$, we have
\bea
\label{EQ13}
		\bar{\hP}^i\circ(\bar{\xi}^i)^{-1}(E^{A_1,\dots, A_n}_{t_1,\dots,t_n})
		\neg=\neg\bar{\hP}^i \big\{ \bigcup_{m=1}^\infty \bigcap_{j=1}^n E^i_{j,m}\big\}\neg=\neg\bar{\hP}^i \big\{ \bigcup_{m=1}^\infty 
		\hat{E}^i_{m}\big\}\neg=\neg\lim_{m\to \infty} 	\bar{\hP}^i \big\{\hat E^i_{m}\big\}\neg=\neg\lim_{m\to \infty} \bar{\hP}^{i,(m)}\big\{\hat E^i_{m}\big\}.
	\eea
Now, by Lemma \ref{alpha}, for two Markovian weak solutions satisfying $\bar{\hP}^i\{\bar{\a}^i_t=\a_t, ~ t\in[0,T)\}=1$, $i=1,2$, we  must have $\bar{\a}^1_t=\bar{\a}^2_t=\a_t=u(t, \Phi(t, \cd))$, $t\in[0,T)$, $\bar{\hP}^1$, $\bar{\hP}^2$-a.s. for some functions $u\in \hL^0([0,T]\times\hR^2)$ and $\Phi\in\hL^0_{\sA}([0,T]\times\O)$. In other words,   $(\bar{\hP}^{i, (m)}, \bar{\xi}^{i,(m)})$, $i=1,2$, satisfy the same SDE (\ref{SDEVXpbar1})  on $[0, \t_m]$ with the same coefficients induced by a (bounded) process $\a^{(m)}$, for which the pathwise uniqueness holds. We conclude that
$\bar{\hP}^{1, (m)}\circ(\bar{\xi}^{1,(m)})^{-1}=\bar{\hP}^{2, (m)}\circ(\bar{\xi}^{2,(m)})^{-1}$.  Note that $\{\t_m^i \ge t_j\} = \{u(t_j, \bar{\xi}^{i,(m)}) \le m\}$, we see that 
$\bar{\hP}^{1,(m)}\big\{\hat E^1_{m}\}=\bar{\hP}^{2,(m)}\big\{\hat E^2_{m}\}$, $m\in \hN$, and the result follows from (\ref{EQ13}).
\qed

\section{Affine Structure of Insider Strategy}
\setcounter{equation}{0}

In the rest of the paper we shall use the STPBVP to construct the equilibrium strategy. Note that 
the solution to STPBVP (\ref{TPBVP}) depends on the ``pricing rule"  $(\m, \rho)$, we first argue that $(\m, \rho)$ can be chosen so that  the equilibrium strategy takes a particular form. 
Specifically,  from Propositions \ref{prop37} and \ref{SDEVXpbar} we see that the $\a$-component in 
a weak solution is closely related to  an ill-posed parabolic PDE (\ref{pdephi}), and in light of the well-known Widder's Theorem and its extensions (cf. e.g., \cite{BPSV, G62, W44, W53}), we may assume that 
$\varphi(t, v, x)=\exp\{ I(t, v,x)\}$, where $I(t, \cd, \cd)$ is quadratic in $(v,x)$. Thus, if a Markovian strategy $\bar{\a}_t=u(t,\Phi(t, \cd))$ (see Remark \ref{alpha}), then  
 \vspace{-2mm}
 \bea
 \label{linearu}
 u(t,v,x) = \rho(t,x)\pa_x \ln\varphi=u_0(t,x)+u_1(t,x)v, \qq (t, v,x)\in[0,T)\times \hR^2,
 \eea
for some functions $u_0,u_1: [0,T]\times \hR \to \hR$ to be determined later. In what follows we call a function $u$ of the form (\ref{linearu}) as having an {\it Affine Structure}. 

We should note that the affine structure of the insider strategy has been widely observed in the literature. 
In particular, the equilibrium strategy  of the form 
\bea
\label{linearstra}
\a_t=\beta_t(V_t-P_t), \qq t\in[0,T), 
\eea
where $\beta=\{\beta_t\}$ is a deterministic function known as the ``trading intensity", can be found in many static information case
(see, e.g., \cite{ABO,  K85}), as well as
dynamic information case (see, e.g., \cite{MSZ}).  The general form in \eqref{linearu} can also be found in \cite{KB92, BP98}. 
In order to validate the affine structure, let us begin with some simple analysis. 

Assume, for example, that a solution to the STPBVP (\ref{SDEVXpbar}) is such that 
$\bar{\a}_t=u(t, \bar{V}_t, \bar{X}_t)$, where $u(t, v,x)$ satisfies (\ref{linearu}), then the function $\varphi$ must have the form $\varphi(t,v,x) = \exp\{I(t,v,x)\}$, where
\bea
\label{EQ17}
I(t,v,x) = h(t,v)+ A(t,x)+ B(t,x)v,
\eea
and $A(t,x)$ and $B(t,x)$ are defined respectively by 
\bea
\label{AB}
A(t,x):=\int_0^x \dfrac{u_0(t,y)}{\rho(t,y)} dy;\q B(t,x):=\int_0^x \dfrac{u_1(t,y)}{\rho(t,y)} dy, \q h(t,v) := \ln\varphi(t,0,v).
\eea
Now assume that $\varphi$ satisfies the PDE (\ref{pdephi}), then we can derive a PDE for the function $I$:
\bea
\label{pdeI}
\left\{\ba{lll}
        I_t + b(t,v,x)I_v +\mu(t,x) I_x+\dfrac{1}{2}\sigma^{2}(t,v,x)[(I_v)^2+I_{vv}]+\dfrac{1}{2}\rho^2(t,x)[(I_x)^2 + I_{xx}]=0;\ms\\
        I(0,v,x)=h(0,v)+A(0,x)+B(0,x)v.
        \ea\right.
\eea
Plugging \eqref{EQ17} into (\ref{pdeI})   we obtain
\bea
\label{EQ18}
  0&=& \dfrac{1}{2}\rho^2(t,x)B_x^2 v^2+ \big\{B_t +\mu(t,x)B_x +\dfrac{1}{2}\rho^2(t,x)[B_{xx} + A_x B_x]\big\}v+A_t +\mu(t,x) A_x\nonumber\\
&& +\dfrac{1}{2}\rho^2(t,x)[A_{xx} +A^2_x]+ h_t + b(t,v,x)[h_v + B ]+\dfrac{1}{2}\sigma^2(t,v,x)\{ h_{vv} + [h_v +B]^2\}.
\eea
For notational simplicity, for given coefficients $b, \si, \m, \rho$, we define
\bea
\label{I012}
\left\{\ba{lll}
I_0(t,x)=I_0(t,x;\m, \rho)=A_t  +\mu(t,x) A_x  +\dfrac{1}{2}\rho^2(t,x)[A_{xx} +A^2_x];\ms\\
I_1(t,x)=I_1(t,x; \m, \rho)=B_t  +\mu(t,x)B_x+\dfrac{1}{2}\rho^2(t,x)[B_{xx} + A_x B_x];\ms\\
I_2(t,x)=I_2(t,x;\m,\rho)= \dfrac{1}{2}\rho^2(t,x)B_x^2;\ms\\
G(t,v,x)= h_t(t,v)+ b(t,v,x)[h_v(t,v)+ B]+\dfrac{1}{2}\sigma^2(t,v,x)\{ h_{vv}(t,v)+ [h_v(t,v)+B]^2\}.
\ea\right.
\eea
Then, (\ref{EQ18}) becomes
\bea
\label{comp}
       I_2(t,x) v^2+I_1(t,x)v+I_0(t,x)+G(t,v,x)=0, \qq (t,v,x)\in[0,T]\times\hR^2.
\eea
We thus obtained the following result for affine structure of function $u$.
\begin{prop}
\label{affine}
The function $u(t,v,x)=\rho(t,x)\pa_x\ln\varphi(t,v,x)$ has an affine structure (\ref{linearu}), where $\varphi$ solves (\ref{pdephi}), if and only if the coefficients $b, \si, \m, \rho$ satisfy the  compatibility conditions (\ref{comp}) with $I_0$-$I_2$ and $G$ being  defined respectively by (\ref{I012}). 

Furthermore, it holds that $\pa_{vvv}G(t, v,x)\equiv 0$, $(t,v,x)\in[0,T]\times\hR^2$.
\qed
\end{prop}

We should note that the compatibility condition (\ref{comp}) is technically difficult to verify in general, as it involves not only a fairly complicated systems of differential equations, but also the selection of the ``pricing rule" $(\m, \rho)$.
In what follows we impose some specific structures on the functions $h$, $b $ and $\si $, and try to find the conditions under which the function $u(t,v,x)$ is of an affine structure. 

Let us begin with an example  of a Kyle-Back problem with dynamic information that fits the generality considered in this paper, and   justifies the validity of the compatibility condition. 
\begin{eg}
\label{egMSZ}
{\rm 
Consider the Kyle-Back problem studied in \cite{MSZ}. More precisely, we assume that
\beaa
\label{eq29}
b(t,v,x) = f_t v + g_t x + k_t, \qq \sigma(t,v,x) = 1.
\eeaa
Denote now $X_t = P_t=\hE^\hP[V_t|\cF^Y_t]$. 
Then, by \cite[Theorem 3.6]{MSZ}, we have
\beaa
\label{eq27}
\mu(t,x) = (f_t+g_t)x + k_t,\qq \rho(t,x)=\rho(t)= S_t\beta_t,
\eeaa
where $S_t$ satisfies a (deterministic) Riccati equation. Furthermore, in \cite{MSZ} it was shown that the equilibrium strategy
 takes the form (\ref{linearstra}).
That is, the equilibrium $\a$ has an affine structure (\ref{linearu}) 
with $u_0(t,x) = -\beta_t x$, $u_1(t,x) = \beta_t$. By definition (\ref{AB}) we then have
\beaa
\label{AB1}
\left\{\ba{lll}
\dis A(t,x) = \int_0^x \frac{u_0(t, y)}{\rho(t,y)}dy=-\frac{1}{S_t}\int_0^x ydy=-\dfrac{x^2}{2S_t};\ms\\
\dis B(t,x) = \int_0^x\frac{u_1(t,y)}{\rho(t,y)}dy=\int_0^x\frac1{S_t}dy=\dfrac{x}{S_t}.
\ea\right.
\eeaa
Plugging these into \eqref{I012} and using the fact that $S$ satisfies the Riccati equation 
$ \frac{dS_t}{dt} = 2f_tS_t -\beta^2_tS^2_t+1$, $ t\in[0,T)$,
one can check that the compatibility condition (\ref{comp}) holds.
\qed}
\end{eg}

In the general nonlinear case, the analysis becomes too complicated to have a generic result. We therefore consider several special cases that might be useful in practice. 

\ms
\no{\bf Case 1.}  $h = h(t)$, $b(t,x,v) = b(t,x)$ and $\sigma(t,v,x)=\sigma(t,x)$. 
 In this case, \eqref{EQ18} is reduced to 
 \begin{align}
 	\label{eq9}
 	I_0 (t,x)+ I_1(t,x )v + I_2(t,x) v^2 = 0,
 \end{align}
where (suppressing variables)
\beaa
	I_0 &=& \pa_t h + bB+ \dfrac{1}{2}\sigma^2B^2+\pa_t A+\mu\pa_xA+\dfrac{1}{2}\rho^2(\pa_{xx}A +(\pa_x A)^2),\\
	I_1 &=& \pa_t B + \mu \pa_x B + \dfrac{1}{2}\rho^2(\pa_{xx}B+\pa_xA \pa_xB),\q I_2  = \dfrac{1}{2}\rho^2 (\pa_x B)^2.
\eeaa
Clearly, \eqref{eq9} implies that $I_0 = I_1=I_2 = 0$. Then, by definition we have $\pa_x B = \frac{u_1(t,x)}{\rho(t,x)}=0$, which implies
that $u_1(t,x) \equiv 0$. Consequently, $B(t,x) \equiv 0$. It then follows that 
 \be
 \label{eq13}
  \pa_t h + \pa_t A+\mu\pa_xA+\dfrac{1}{2}\rho^2(\pa_{xx}A +(\pa_x A)^2)= 0.
 \ee
 That is, a necessary condition for  affine structure is that $u_1 \equiv 0$ and   $h$, $u_0$, $b, \sigma, \mu, \rho$  satisfy \eqref{eq13}.
 
\ms

\no{\bf Case 2.}   $h = h(t)$, $b(t,v,x) = b_0(t,x)+b_1(t,x)v$, $\sigma(t,v,x) = \sigma_0(t,x)+ \sigma_1(t,x)v$. Then, similar to Case 1, 
we simplify the equation \eqref{EQ18} and denote the coefficients as $I_0, I_1,I_2$, where 
\beaa
\left\{\ba{lll}
	I_0 &=& \pa_t h(t)+ b_0B+ \dfrac{1}{2}\sigma_0^2B^2+\pa_t A+\mu\pa_xA+\dfrac{1}{2}\rho^2(\pa_{xx}A +(\pa_x A)^2);\\
	I_1 &=& b_1B+ \sigma_0\sigma_1B^2+ \pa_t B + \mu \pa_x B + \dfrac{1}{2}\rho^2(\pa_{xx}B+\pa_xA \pa_xB);\\
	I_2 &=& \dfrac{1}{2}\rho^2 (\pa_x B)^2 +  \dfrac{1}{2}\sigma_1^2 B^2.
	\ea\right.
\eeaa
We see from $I_2 =0$ that $u_1 \equiv 0$, which again leads to \eqref{eq13}.

\ms

\no{\bf Case 3.}  $h = h(t)$, $b(t,v,x) = b_0(t,x) + b_1(t,x)v+b_2(t,x)v^2$, $\sigma(t,v,x) = \sigma_0(t,x)+ \sigma_1(t,x)v$. 
In this case we need 
\bea
\label{eq18}
\left\{\ba{lll}
I_0 =\pa_t h(t) + b_0B+\dfrac{1}{2}\sigma_0^2B^2+\pa_t A+\mu\pa_xA+\dfrac{1}{2}\rho^2(\pa_{xx}A +(\pa_x A)^2)=0;\\
I_1 =b_1B + \sigma_0\sigma_1B^2 + \pa_t B + \mu \pa_x B + \dfrac{1}{2}\rho^2(\pa_{xx}B+\pa_xA \pa_xB)=0;\\
I_2 =\dfrac{1}{2}\rho^2 (\pa_x B)^2 +  \dfrac{1}{2}\sigma_1^2 B^2 + b_2B =0.
\ea\right.
\eea
In particular, $I_2=0$  if and only if
\be
\label{eq15}
u_1^2(t,x) = -\sigma_1^2\Big( \int_{x_0}^x\dfrac{u_1(t,y)}{\rho(t,y)}dy\Big)^2 - 2b_2 \int_{x_0}^x\dfrac{u_1(t,y)}{\rho(t,y)}dy.
\ee

 If we choose $u_1  = \rho $, then \eqref{eq15} implies
$\rho^2 = -\sigma^2_1 (x-x_0)^2-2b_2(x-x_0)$. Using $I_1=0$ in \eqref{eq18}, we can write $u_0$ as 
\bea
\label{eq17}
u_0 = \dfrac{2}{u_1}\big[-\pa_t B -\mu \pa_xB -\dfrac{1}{2}\rho^2\pa_{xx} B -b_1B - \sigma_0\sigma_1B^2 \big].
\eea
Therefore, \eqref{eq18}, together with \eqref{eq15}, \eqref{eq17}, guarantees the affine structure in this case.

\ms
\no{\bf Case 4.} $h(t,v)=h_0(t) + h_1(t)v$,  $b, \sigma$ same as Case 3. 
In this case,
\beaa
\label{eq20}
\left\{\ba{lll}
	I_0 &=& \pa_t h_0  + b_0(h_1 + B) +\dfrac{1}{2}\sigma_0^2(h_1 +B)^2+\pa_t A+\mu\pa_xA+\dfrac{1}{2}\rho^2(\pa_{xx}A +(\pa_x A)^2)=0;\\
	I_1 &=&  \pa_t h_1 +  b_1(h_1  +B) + \sigma_0\sigma_1(h_1 + B)^2 + \pa_t B + \mu \pa_x B + \dfrac{1}{2}\rho^2(\pa_{xx}B+\pa_xA \pa_xB)=0;\\
	I_2 &= & \dfrac{1}{2}\rho^2 (\pa_x B)^2 + \dfrac{1}{2} \sigma_1^2 (h_1 +B)^2 + b_2(h_1  +B) =0.
\ea\right.
\eeaa

\ms
\no{\bf Case 5.} $h = h_0(t) + h_1(t)v + h_2(t)v^2$. 
Since there are the terms $b h_v$, $\sigma^2 \pa_vh^2$ in $G(t,v,x)$, and $h$ is quadratic, we conclude that $\sigma(t,v,x)$ must be independent of $v$, and $b$ has to be linear in $v$.
We thus assume that $b=b_0(t,x)+ b_1(t,x)v$, $\sigma = \sigma(t,x)$, in other words, 
\beaa
	\label{eq22}
\left\{\ba{lll}
	I_0 = \pa_t h_0  + b_0(h_1 + B) +\dfrac{1}{2}\sigma^2[2h_2 +(h_1 +B)^2]+\pa_t A+\mu\pa_xA+\dfrac{1}{2}\rho^2(\pa_{xx}A +(\pa_x A)^2)=0;\\
	I_1 =  \pa_t h_1 + 2b_0h_2 +  b_1(h_1  +B) + 2\sigma^2(h_1 + B)h_2 + \pa_t B + \mu \pa_x B + \dfrac{1}{2}\rho^2(\pa_{xx}B+\pa_xA \pa_xB)=0;\\
	I_2 = \pa_t h_2 + 2 b_1 h_2 +2\sigma^2 h_2^2+\dfrac{1}{2}\rho^2 (\pa_x B)^2  =0.
	\ea\right.	
\eeaa

\section{The Filtering Problem and FBSDE under Affine Structure}
\setcounter{equation}{0}

A popular approach in studying Kyle-Back equilibrium problem is nonlinear filtering (cf. e.g., \cite{ABO, Dani, MSZ}). In fact,  in the dynamic information case where the market price is in the form of an {\it optional projection}:  $P_t = \hE[V_t | \cF^Y_t]$, $t\in[0,T]$, we believe that the filtering approach should be particularly effective in determining the equilibrium strategy, which we now explain. 
 
We begin by recasting the STPBVP (\ref{TPBVP}) as a nonlinear filtering problem. Consider a (Markovian) weak solution $(\bar{\hP}, V, X, B, \a)$ to STPBVP (\ref{TPBVP}), where $ \a _t=u(t,  {V}_t,  {X}_t)$, and under $\bar\hP$, 
	\bea
	\label{filtereq}
		\left\{\begin{array}{lll}
		dV_t = b(t, V_t, X_t)dt +\sigma(t, V_t, X_t) dB^1_t, \qq \qq\qq\qq& V_0=v_0; \ms\\
		dX_t = [\mu(t, X_t)+\rho(t,X_t)u(t, V_t, X_t)]dt +\rho(t,X_t)dB^2_t, & X_0=x_0; \ms\\
		dY_t=u(t, V_t, X_t)dt+dB^2_t, & Y_0=0. 
		\end{array}\right.
	\eea
Since the function $u$ is now fixed, (\ref{filtereq}) can be thought of as a nonlinear filtering problem with  correlated noises, in which $(V,X)$ is the signal process and $Y$ is the observation process. The only technical problem, however, is whether the function $u$ satisfies usual technical requirements so that the SDE for $P_t=\hE[V_t|\cF^Y_t]$ (known as the 
Fujisaki-Kallianpur-Kunita (FKK) equation \cite{FKK}) holds. To this end, we assume 
that $u(t,v,x)$ has the {\it affine structure}:
 $u(t,v,x) = u_0(t,x)+u_1(t,x)v$. Denoting $\a_t=u(t, V_t, X_t)$, and consider the  SDE:
 \bea
 \label{M}
		dM_t = -\alpha_t M_t dB^2_t,\quad M_0 = 1, \qq t\in[0,T].
\eea
The following result is a modification of 
 \cite[Lemma 4.1.1]{Ben} to the current case.
\begin{prop}
\label{Mmg}
	Assume Assumptions \ref{assump1}, and assume further that the function $u$ in 
(\ref{filtereq}) satisfies  $|u(t, v,x)| \le K(t)(1 +|v|+|x|)$, for $(t,v,x)\in [0,T)\times \hR^2$, for some function $K \in \hL^2([0,T];\hR_+)$. 
Then, the solution $M$ to (\ref{M}) is a true martingale on $[0,T]$. 
\end{prop}
{\it Proof.} Clearly, $M$ is a local martingale. 
 Then,  by Fatou's lemma, for any $t\in [0,T]$, we have $\hE[M_t] \leq\lim_{n\to \infty} \hE[M_{t \wedge \tau_n}]=\hE[M_0]=1$, where $\{\t_n\}$ is any announcing sequence for $M$, and $M$ is a true martingale iff 
$\hE[M_t] = 1, 0\leq t \leq T$, which we now prove. 
For any $\e> 0$, define  $f_{\e}:= \frac{x}{1+\e x}$, and $M^\e_t:=f_\e(M_t)$, $t\in[0,T]$. Clearly, by bounded convergence theorem, we have $\lim_{\e\to0}\hE[M^\e_t]  = \hE[M_t]$. On the other hand, by a simple application of It\^o's formula and then taking expectation one has 
	\bea
	\label{Mepsi}
		\hE[M^\e_t]
		:=\dfrac{1}{1+\e}-\hE \Big[ \int_0^t G^\e(\a_s, M_s)ds\Big], ~t\in[0,T],
	\eea
where $G^\e(\a, x):=\frac{\e \a^2 x^2}{(1+\e x)^3}$. It is easy to check that there exists $C>0$, such that $|G^\e(\a, x)|\le C \a ^2  x $, for all $\e, x>0$. Denoting $U_t := M_t(V_t^2 + X_t^2)$, then the linear growth assumption for $\a_t$ gives
\begin{center}
	$\hE[G^\e (\a_t, M_t)] \le C \hE [\a_t^2 M_t] \le CK^2(t) \big[1 + \hE[U_t]\big].$
\end{center}
We  claim that
$\sup_{t\in[0,T]}\hE [ U_t] <\infty$. The result then follows easily from the Dominated Convergence theorem. Applying It\^o's formula to $U_t$ and $f_\e(U_t)$, we have (denoting $|\xi_0|^2=v_0^2+x_0^2$)
	\begin{align*}
		f_\e(U_t)
		&= \dfrac{|\xi_0|^2}{1+\e|\xi_0|^2}+\int_0^t\dfrac{2M_s\big[V_sb_s+X_s\mu_s+\frac{1}{2}(\sigma_s^2+\rho_s^2)\big]}{(1+\e U_s)^2}ds+\int_0^t \dfrac{2M_sV_s\sigma_s }{(1+\e U_s)^2}dB_t^1\\
		&+\int_0^t \dfrac{-\e\big[4V_s^2\sigma_s^2M_s^2+\big(2M_sX_s\rho_s-U_s\alpha_s\big)^2\big]}{(1+\e U_s)^3}ds +\int_0^t \dfrac{-U_s\alpha_s + 2M_sX_s\rho_s}{(1+\e U_s)^2} dB^2_s.
	\end{align*}
Taking expectation on both sides, and by the linear growth of $b,\sigma, \mu$ and $\rho$, we obtain
\begin{align*}
\hE [f_\e (U_t) ] 
		\leq |\xi_0|^2+\int_0^t\hE\Big[\dfrac{2M_s\big[V_sb_s+X_s\mu_s+\dfrac{1}{2}(\sigma_s^2+\rho_s^2)\big]}{(1+\e U_s)^2}\Big]ds
		\leq |\xi_0|^2+ \int_0^t L ( \hE [f_\e (U_t)] +1 )ds.
\end{align*}
Now, first applying Gronwall's inequality  and then applying Fatou's lemma (sending $\e\to0$), we deduce that $\sup_{t\in[0,T]}\hE[U_t ]<\infty$, proving the claim.
\qed

  We should note that with Proposition \ref{Mmg} and the affine structure assumption on $u$ the SDE (\ref{filtereq}) can be
naturally extended to $[0,T]$, and we can follow the same argument of \cite[Theorem 4.1]{FKK}
to derive the FKK equation for $P_t=\hE^{\bar\hP}[V_t|\cF^Y_t]$, which takes the following form:
\bea
\label{P_SDE}
\left\{\ba{lll}
dP_t =[\hE^t[b(t,V_t,X_t)]- \hE^t[u(t,V_t,X_t)] Z_t ]dt + Z_tdY_t, \ms\\
	Z_t :=\hE^t[V_t u(t,V_t,X_t)] - P_t \hE^t[u(t,V_t,X_t)], 
\ea\right. \qq t\in[0, T],
\eea
where $\hE^t[\cd]:=\hE^{\bar\hP}[\cd|\cF^Y_t]$, $t\in[0,T]$. Now if we assume  that the coefficient $b(\cds)$ is also of affine structure:
$b(t,v,x) = b_0(t,x)+b_1(t,x)v$, and $X$ is $\hF^Y$-adapted, then \eqref{P_SDE} can be rewritten as
\bea
\label{lp}
	dP_t = \{b_0(t,X_t)+b_1(t,X_t)P_t-(u_0(t,X_t)+u_1(t,X_t)P_t ) Z_t\}dt + Z_t dY_t, \qq t\in[0,T]. 
\eea

Let us now choose $\a_t=u(t, V_t, X_t)$, $t\in[0,T]$, to be the $\a$-component of a Markovian weak solution to the STPBVP (\ref{TPBVP}), and assume that it has the affine structure. 
By Proposition \ref{Mmg}, the process $M$ defined by (\ref{M}) is a martingale on $[0,T]$, so we can define a new probability measure $\bar\hQ$ on the canonical space $(\O, \cF)$ by $\frac{d\bar{\hQ}}{d\bar{\hP}}|_{\cF_T}=M_T$, then under $\bar{\hQ}$, the process
$Y$ (for the given $\a$) is a Brownian motion, and $\bar{\hQ}\{V_T=g(X_T)\}=\bar{\hP}\{V_T=g(X_T)\}=1$. In other words, under 
$\bar\hQ$, we can rewrite \eqref{lp} and the SDE (\ref{TPBVP}) for $X $ as the following forward-backward SDE (FBSDE): 
	\begin{equation}
	\label{XP_SDE}
	\left\{\begin{array}{lll}
	dX_t = \mu(t,X_t) dt + \rho(t,X_t) dY_t, \q X_0=x;\ms\\
		dP_t = [\beta_0(t,X_t, P_t)+\beta_1(t,X_t, P_t ) Z_t]dt + Z_t dY_t, \q 	P_T = g(X_T), 
	\end{array}\right.
	\end{equation}
where $ \beta_0(t,x,y)=-b_0(t,x) - b_1(t,x)y$, $\beta_1(t,x, y) =u_0(t,x)+ u_1(t,x)y$.
\begin{rem}
\label{remark4.0}
{\rm
(i) Although $\bar{\hQ}\sim \bar{\hP}<\neg\neg<\hQ^0$ and the process  $Y$ is a Brownian motion under both measures $\bar\hQ$ and $\hQ^0$, $\bar\hQ$ and $\hQ^0$ are not equivalent on $\cF_T$, since $\hQ^0\{V_T\neq g(X_T)\}>0$ in general.  In fact, $L^\nu$ is local martingale, but $M$ is a true martingale. 

\ss
(ii) Under Assumption \ref{assump1},  $X$ is a diffusion driven by the $\bar\hQ$-Brownian motion $Y$, hence it is $\hF^Y$-adapted, which justifies (\ref{lp}), whence (\ref{XP_SDE}). 
\qed}
\end{rem} 

We should note that the FBSDE \eqref{XP_SDE} is actually ``decoupled", in the sense that the forward SDE is independent of the backward components $(Y, Z)$. But the BSDE in (\ref{XP_SDE}) is somewhat non-standard in that the coefficients are neither Lipschitz nor of linear growth. Specifically, the fact that $|\beta_1(t,x,y)z|\le K(1+|y||z|)$ makes it super-linear in $(y, z)$, and is beyond the usual ``quadratic BSDE" framework. Nevertheless, the well-posedness of \eqref{XP_SDE} can be argued via a 
more or less standard localization argument following the idea of \cite{MaLiu}. Since this is not the main purpose of the paper, 
we shall only state the following result, but omit the proof (see \cite{TanY} for details). 
\begin{prop}
\label{prop4.2}
Assume Assumption \ref{assump1}, and let $(\bar{\hP}, (B^1, B^2), (V, X), \a)$ be a Markovian nested solution to STPBVP (\ref{TPBVP}), and assume that $\a$ has an affine structure. Then there exists a probability measure $\bar\hQ$ on the canonical space $(\O, \cF)$,  such that

(i) $\dfrac{d\bar\hQ}{d\bar\hP}\big|_{\cF_T}=M_T$, where $M$ satisfies the linear SDE (\ref{M});

\ss
(ii) denoting $Y_t=B^2_t+\int_0^t \a_sds$ and $P_t=\hE^{\bar\hP}[V_t|\cF^Y_t]$, $t\in[0,T]$, then $Y$ is a $\bar{\hQ}$-Brownian motion, and under $\bar\hQ$, $(X, P)$ satisfies the FBSDE (\ref{XP_SDE}). 
\qed
\end{prop}

In the rest of this section we  try to determine the most important element of the pricing mechanism: the function $H:[0,T]\times \hR\mapsto \hR$, so that $P_t=H(t, X_t)$, $t\in[0,T]$. 
To begin with, we recall from the general theory of FBSDE (cf. e.g., \cite{MYbook}, \cite{MWZZ})  that, if $(X, P, Z)$ is the solution to the FBSDE (\ref{XP_SDE}), then under appropriate conditions on the coefficients, there is a  {\it decoupling field} $H:[0,T]\times\hR\mapsto \hR$, which satisfies the following semilinear PDE (at least in the viscosity sense):
\begin{equation}
\label{pdeh}
\left\{\begin{array}{lll}
H_t(t,x)+ \dfrac{1}{2}\rho^2(t,x)H_{xx}(t,x) + \m(t,x)H_x(t,x) +h(t, x, H(t,x), \rho(t, x)H_x(t, x))=0;\\
H(T,x)=g(x), 
\end{array}
\right.
\end{equation}
where $h(t, x, y, z)=\beta_0(t, x, y)+\beta_1(t, x, y)z$,  and the following relation holds:
$P_t=H(t, X_t)$, $Z_t=\rho(t, X_t)H_x(t, X_t)$, $t\in[0,T]$.
The following extension of Example \ref{egMSZ} justifies this fact.
\begin{eg}
	\label{eg2}
	{\rm
		Recall that in Example \ref{egMSZ}, in which the coefficients $b, \sigma, \mu$ and the function $u$ have the specific form: $b(t,v,x)=f_tv+g_tx+k_t$, $\si\equiv 1$, $\mu(t,x) = (f_t+g_t)x+k_t$, $u(t,v,x)=\beta_tv-\beta_t x$, and thus the PDE (\ref{pdeh}) now reads
		(suppressing variables): 
		\bea
		\left\{\begin{array}{lll}
			\label{eq601}
			H_t +\big((f_t+g_t)x+k_t+\rho(-\beta_tx+\beta_tH )\big)H_x +\dfrac{1}{2}\rho^2H_{xx}  = g_tx+k_t+f_tH ;\ms\\
			H(T,x)=x,
		\end{array}\right.	
		\eea
		We can easily check that $H(T,x) =x$ is the (unique) solution to (\ref{eq601}), and hence $P_t =H(t,X_t)=X_t$, for $t\in[0,T)$, and 
		$X_T =H(T, X_T)=P_T=V_T$. 
		\qed}
\end{eg}

\begin{rem}
\label{remark6.1}
{\rm 	If we restrict the strategy to the form $\a _t = \beta_t (V_t - P_t) = \beta_t(V_t-H(t,X_t))$,  that is, $u_0 = - \beta_t H(t,x), u_1 = \beta_t$, and we assume further that the original asset $V$ is under the risk neutral probability so that  $b(t,v,x)=0$, then \eqref{pdeh} is reduced to 
	\bea
	\label{Hpde0}
	\left\{\begin{array}{lll}
		H_t(t,x)+\mu(t,x) H_x(t,x)+\dfrac{1}{2}\rho^2(t,x)H_{xx}(t,x) = 0;\ms\\
		H(T,x)=g(x).
	\end{array}\right.	
	\eea
We should note that the PDE (\ref{Hpde0}) is well-posed with properly given coefficients, and $(\m, \rho)$ can be chosen as part of the pricing rule. The determination of $(\m, \rho)$, however, is the main task for finding the Kyle-Back equilibrium, which will be discussed in details in the next section. 
\qed
}
\end{rem}

\section{Sufficient Conditions for Optimality}
\setcounter{equation}{0}

We are now ready to investigate the main issue of this paper: finding the equilibrium of the pricing problem. That is, we are to find the optimal strategy $\alpha^*$ for the insider, which maximizes her expected terminal wealth $W_T$, given the  pricing rule $P_t=\hE[V_t|\cF_t]$, $t\in [0,T]$. 

In light of the analysis in the previous sections, we can recast the problem of finding the Kyle-Back equilibrium as follows. First recall the Markovized system (\ref{vx00}):
\bea
\label{vxy}
	\left\{\begin{array}{lll}
	dV_t = b(t, V_t, X_t)dt + \sigma(t, V_t, X_t) dB^1_t, \qq & V_0=v; \ms\\
	dX_t =\mu(t, X_t)dt+\rho(t,X_t)dY_t= [\mu(t, X_t)+\rho(t,X_t)\a_t]dt +\rho(t,X_t)dB^2_t, & X_0=x.
	\end{array}\right.
\eea
where $\a\in \sU_{ad}$ (see (\ref{Uad}) for definition). Assume that the process $\a$ takes the feedback form $\a_t=u(t, V_t, X_t)$, we have argued in \S2 that finding the optimal strategy amounts to solving a stochastic control problem with state equation (\ref{vxy}) (or (\ref{vx00})) and the cost functional (\ref{J0}). Moreover, a necessary condition for $\a \in \sU_{ad}$ being an equilibrium is that 
$V_T=P_T=H(T, X_T)=g(X_T)$ (see (\ref{terminal})). Therefore, 
We shall  consider only  the (weak) solution $( \bar\hP,  V , X , \a)$ to  STPBVP (\ref{SDEVX}), and by 
Proposition \ref{limitbarP}, we shall assume that 
 $\bar{\hP}|_{\cF_t}<\neg\neg<\hQ^0|_{\cF_t}$, $t<T$, and $\bar{\hP}\circ (V_{T})^{-1} =m^*$.
 
It is worth noting that  the candidates for the equilibrium described above is based on the SDE (\ref{vxy}), 
whence 
the coefficients $(\m, \rho)$. We shall argue that the equilibrium can be determined by properly choosing $(\m, \rho)$ through some ``compatibility conditions".
 
 \ms
 
 \no{\bf The case $b(t,v,x) \equiv 0$.} For notational simplicity, in what follows we use $\hP$ instead of $\bar{\hP}$. As we pointed out in Remark \ref{remark6.1}, this could be the case when $\hP$ is the risk neutral probability measure, and $V$ is the discounted asset price, hence a $(\hP, \hF)$-martingale. We note that in this case the market price $P_t=\hE[V_t|\cF^Y_t]$,  $t\ge 0$ is a $(\hP, \hF^Y)$-martingale. Indeed, since $V=\{V_t\}$ is a $(\hP, \hF)$-martingale, for $s<t$, we have
	\beaa
	P_s = \hE^\hP[V_s | \cF^Y_s] = \hE^\hP[ \hE^\hP[V_t|\cF_s]|\cF^Y_s] = \hE^\hP[V_t | \scF^Y_s]= \hE^\hP[ \hE^\hP[V_t|\cF^Y_t]|\cF^Y_s]=\hE^\hP[P_t | \cF^Y_s].
	\eeaa
On the other hand, if we assume that $P_t=H(t, X_t)$, $t\in[0,T]$, where $X$ satisfies the SDE (\ref{vxy}), such that $P_T=g(X_T)$, then a simple application of It\^o's formula shows that $P=\{P_t\}$ begin an $\hF^Y$-martingale means that the decoupling field $H$ must satisfy the PDE:
\bea
\label{eq55}
\left\{\begin{array}{lll}
	H_t+ \mu(t,x)\pa_xH+\dfrac{1}{2}\rho^2(t,x)H_{xx}= 0;\qq t\in[0,T)\ms\\
	H(T,x)=g(x).
\end{array}\right.	
\eea
Comparing to (\ref{pdeh}) and recalling (\ref{XP_SDE}) we see that under the affine structure we should have 
$$ h(t, x, H(t, x), \rho(t,x)H_x(t, x))=\beta_1(t, x, H)\rho(t,x)H_x=(u_0(t,x)+ u_1(t,x)H)\rho(t,x)H_x\equiv 0.
$$ 
Therefore we have $u_0(t,X_t) = - \beta_t H(t,X_t)$, where $ \beta_t =u_1(t,X_t)$. Consequently, we see that $\a_t=u_0(t, X_t)+u_1(t, X_t)V_t=u_1(t, X_t)(V_t-H(t, X_t))=\beta_t(V_t-P_t)$, which is exactly the form commonly seen in the literature (see, e.g., \cite{ABO,K85,MSZ}),  except that $\beta_t$ is no longer deterministic.
Our first main result of this section is the following theorem. 
\begin{thm}
	\label{thm1}
	Assume Assumption \ref{assump1}, and that $b(t,v,x)=0$. Assume further that $(\bar\hP,  \bar{V}, \bar{X}, \bar{\a})$ is a weak solution to the STPBVP (\ref{SDEVX}) such that $\bar\a_t$ has an affine structure. Then, 
	
	(i) the market price $P_t=\hE^{\bar\hP}[\bar{V}_t|\cF^Y_t]=H(t,\bar{X}_t)$, $t\in[0,T)$ is an $\hF^Y$-martingale, where
	$H$ solves the PDE (\ref{eq55});
	
\ss

(ii) the process $\bar\a$ is of the form $\bar{\a}_t= \beta(t,\bar{X}_t)(\bar{V}_t - H(t,\bar{X}_t))=\beta(t, \bar{X}_t)(\bar{V}_t-P_t)$, $t\in[0,T)$, where $(V,X)$ is the solution to the SDE
	(\ref{vxy}) under some probability measure $\bar\hP$, such that  $\bar{V}_T = g(\bar{X}_T)$, $\bar\hP$-a.s.; 

\ss
(iii) $\bar\a$ is an equilibrium strategy if the following ``compatibility condition" holds for $(\m, \rho)$: 
	\begin{align}
		\label{eq57}
		\pa_t\rho(t,x)-\pa_x\mu(t, x)\rho(t, x)+\pa_x\rho(t, x)\mu(t,x)+\dfrac{1}{2}\rho^2(t,x)\pa_{xx}\rho(t, x) = 0.
	\end{align} 
\end{thm}

{\it Proof.} The parts (i) and (ii) have been argued prior to the theorem. We shall prove only part (iii). To this end, we shall borrow the idea of \cite{WU}, and look for a function  $J(t,x;a) $ such that  for fixed $a\in\hR$, $J(\cd, \cd;a)\in \hC^{1,2}([0,T]\times\hR)$, and satisfies the following properties
\bea
\label{PDEJ}
\left\{\ba{lll}
J_t(s,x;a)+J_x(s,x;a)\mu(s,x)+\dfrac{1}{2}J_{xx}(s, x;a)\rho^2(s,x)= 0; \ms\\
J_x(s,x;a)\rho(s,x) = H(t,x)-a; \ms\\
		J(T,x;a) \geq 0, \mathrm{and} \,\, J(T,x;a) = 0\,\, \mathrm{iff}\,\, a = g(x).
\ea\right.
\eea

Assume now that  a function $J$ satisfying (\ref{PDEJ}) exists. Then for any $\a\in\sU_{ad}$, we let $(\hP, V^\a, X^\a)$ be a weak solution to the SDE (\ref{SDEVX}).  
Given $a\in\hR$, applying It\^o's formula to $J(\cd, \cd;a)$  we have
	\bea
	\label{eq12}
		J(t,X^\a_t;a) & =& J(0,x_0;a)+\int_{0}^{t} \big[J_t(\cd, \cd;a)+J_x(\cd, \cd;a)\mu+\dfrac{1}{2}J_{xx}(\cd, \cd;a)\rho^2\big] (s,X^\a_s)ds\nonumber\\
		&&+ \int_{0}^{t} J_x(s,X^\a_s;a)\rho(s,X_s) dY_s=J(0,x_0;a)+\int_{0}^{t} (H(s,X^\a_s) -a )dY_s \\
		&=&J(0,x_0;a)+\int_{0}^{t} (H(s,X^\a_s) -a )\alpha_s ds+\int_{0}^{t} H(s,X^\a_s) dB^2_s - aB^2_t. \nonumber
	\eea
Denoting $(V,X)=(V^\a, X^\a)$ and applying the total probability formula and (\ref{eq12}) we have
	\bea
	\label{eq131}
	&&	\hE^{ \hP}\big[J(T,X_T;V_T) -J(0,x_0;V_T) \big] 
		=\int_\hR \hE^{ \hP}\big[J(T,X_T;a) -J(0,x_0;a) |  V_T=a\big] \hP\circ (V_T)^{-1} (da)\nonumber\\
		& =& \int_\hR  \hE^{\hP}\Big[\int_{0}^{T} (H(s,X_s) - a) \alpha_tdt+ \int_{0}^{T} H(t,X_t) dB^2_t - aB^2_T | a=V_T\Big] \hP\circ (V_T)^{-1}(da)\\
		&=&\hE^{ \hP}\Big[\int_{0}^{T} (H(s,X_s) - V_T) \alpha_tdt\Big]+\hE^{ \hP}\Big[\int_0^T H(t,X_t)dB^2_t\Big]- \hE^{\bar\hP}[V_TB^2_T]\nonumber\\
		&=& \hE^{ \hP}\Big[\int_{0}^{T} (H(s,X_s) - V_T) \alpha_tdt\Big]- \hE^{ \hP}[V_TB^2_T]. \nonumber
			\eea
But, since $\lan B^1, B^2\ran\equiv 0$, we have $d(V_tB^2_t) = V_t dB^2_t +B^2_t \sigma(t, V_t, X_t) dB^1_t$, $t\ge 0$. That is,  $\{V_t B^2_t\}$ is a $\hP$-martingale, hence $\hE^{\hP}[V_TB^2_T] = 0$.
	Recalling (\ref{w2}) we deduce from (\ref{eq131}) that 
	\bea
	\label{EWa}
		\hE^\hP[W^\alpha_T] &=& \hE^{ \hP}\Big[\int_{0}^{T} ( V^\a_T -H(s,X^\a_s)) \alpha_tdt\Big]=\hE^{\hP}\big[J(0,x_0;V^\a_T) -J(T,X^\a_T;V^\a_T) \big]\\
 &\le &\hE^{\hP}\big[J(0,x_0;V^\a_T)\big].\nonumber
 	\eea
Here the last inequality is due to property (\ref{PDEJ}) of the function $J$, and furthermore, the equality holds if and only if the terminal condition $V^\a_T=g(X^\a_T)$ holds. 	Consequently, if we  let $(\bar\hP, \bar{V}, \bar{X}, \bar\a)$ be a weak solution to STPBVP (\ref{SDEVX}), then Proposition \ref{limitbarP}, together with (\ref{EWa}), shows that 
\beaa
\label{EWbara}
\hE^{\bar\hP}[W^{\bar\a}_T]=\sup_{\a\in\sU_{ad},\hP\circ (V^\a_T)^{-1}=m*} \hE^\hP[W^\a_T]=\int_\hR J(0,x_0;a)m^*(da).
\eeaa
In other words, the solution to the STPBVP leads to the optimal strategy for the insider, among all the strategies satisfying  $\hP\circ (V^\a_T)^{-1} = m^*$.
	
	Our last task is to construct a function $J$ that satisfies all the requirements in (\ref{PDEJ}). In light of \cite{WU}, we consider the following function:
\begin{align}
		\label{eq10}
		J(t,x;a) =\int_{g^{-1}(a)}^x \dfrac{H(t,y)-a}{\rho(t,y)}dy +\int_t^T f(s;a)ds,
	\end{align}
	where $H(\cd,\cd)$ satisfies (\ref{eq55}), and $f(t;a)$ is a function to be determined and independent of $x$. To check that such a function is possible for the proper choices of $\m, \rho$, and $f$, we simply plugging the function $J$ into the PDE in (\ref{PDEJ}) to get 
	\begin{align}
		\label{eq11}
		f(t;a) = &\Big[\Big(\dfrac{\mu}{\rho}-\dfrac{\rho_x}{2}\Big)(H-a)\Big](t,x)+\dfrac{(H_x\rho)(t,x)}{2}+\int_{g^{-1}(a)}^x\bigg[\dfrac{H_t}{\rho } - \dfrac{(H-a)\rho_t}{\rho^2 }\bigg](t,y)dy.
	\end{align}
In order that $f(\cd;a)$ is independent of $x$, we take derivative of the right hand side of \eqref{eq11} with respect to $x$, and multiply it by $\rho^2(t,x)$ to obtain (suppressing variables and rearranging terms)
	\beaa
	f_x\rho^2&=& \rho[ H_t+\mu H_x+ \dfrac{1}{2}\rho^2H_{xx}]+[(\mu_x \rho-\mu \rho_x) 
		- \dfrac{1}{2}\rho_{xx} \rho^2  -\rho_t] (H -a)\\
		&=&[(\mu_x \rho-\mu \rho_x) - \dfrac{1}{2}\rho_{xx} \rho^2  -\rho_t] (H -a), 
	\eeaa
thanks to (\ref{eq55}). Since $\rho$ is positive, we see that $f_x\equiv 0$ provided \eqref{eq57} holds. We note that if the function $f$ in
(\ref{eq10}) is independent of $x$, then the second equation in (\ref{PDEJ}) is obvious by definition. It thus remains to verify the last requirement of (\ref{PDEJ}). To see this we note that 
	\begin{equation*}
		\label{eq8}
		J(T,x;a) =\int_{g^{-1}(a)}^{x} \dfrac{H(T, y) - a}{\rho(T,y)} dy= \int_{g^{-1}(a)}^{x} \dfrac{g(y) - a}{\rho(T,y)} dy.
	\end{equation*}
Since $g$ is increasing, for any $y\ge g^{-1}(a)$ we have $g(y)\ge g(g^{-1}(a))=a$. Thus $J(T, x;a)\ge 0$, for $x\ge g^{-1}(a)$ as  $\rho(T, y)>0$,  and $J(T, x; a)=0$ if and only if $x=g^{-1}(a)$. The proof is now complete. 
\qed

\begin{rem}
{\rm The compatibility condition (\ref{eq57}) between the coefficients $\m, \rho$, and the PDE (\ref{eq55}) for the pricing rule $H$ are 
not new. In the so-called ``long-lived" information case, for example, the market price $P_t=\hE[V_T|\cF^Y_t]$, $t\ge 0$, is naturally a martingale, and $b\equiv 0$ is by assumption, thus Theorem \ref{thm1} always applies. In this case, 
\cite{WU} chooses $\m=0$ and $\rho = 1$, which obviously satisfies the compatibility condition (\ref{eq57}), and (\ref{eq55}) becomes 
$H_t+\dfrac{1}{2}H_{xx} = 0$, and $f(t) = H_x(t,g^{-1}(a))$.
	
As another example, in \cite{CCD11} it is derive from a control theoretic argument via HJB equation that $\mu = 0$, and $\rho$ and $H$ satisfy
	\begin{align*}
		\rho_t +\dfrac{\rho^2}{2}\rho_{xx} = 0,\qq
		H_t+\dfrac{\rho^2}{2}H_{xx}=0, 
	\end{align*}
and $f(t;a) = H_x(t,g^{-1}(a))\rho(t,g^{-1}(a))$, justifying our compatibility conditions (\ref{eq55}) and (\ref{eq57}).
	\qed}
\end{rem}

\ms
\no{\bf The Case $b(t, v, x)\neq 0$.}
In this case, the market price $P_t=\hE[V_t|\cF^Y_t]$, $t\ge 0$, is no longer an $\hF^Y$-martingale, but rather an ``optional projection", and the discussion is more involved. Assume again that both $b$ and $\a$ have the affine structure: $b=b_0(t,x)+b_1(t,x)v$ and $u=u_0(t,x)+u(t,x)v$. By Proposition \ref{prop4.2}, the decoupling field $H(t,x)$ would satisfy a more general PDE:
\begin{eqnarray}
	\label{eq58}
		H_t+\mu H_x  +\dfrac{1}{2}\rho^2 H_{xx}= - (b_0+ b_1 H)+ (u_0 + u_1 H)\rho H_x,\qq 
		H(T,x)=g(x).
\end{eqnarray}
So if we still try to construct the function $J(t,x;a)$ as in \eqref{eq11}, then it  may not be possible to find a corresponding function
$f$ that is independent of $x$.  We now try to modify \eqref{eq11} accordingly. First recall \eqref{w1} we can write 
\beaa
\label{EWTa-1}
	\hE^{\hP}[W^\a_T]=\hE^{\hP}\Big[\int^T_0[F(t, V^\a_t,X^\a_t)-H(t, X^\a_t)]u(t,V^\a_t,X^\a_t)dt\Big], 
\eeaa
where $F(t,V_t,X_t) := \hE^{\hP}[ V_T | \cF ^{V,X} _t ]$, thanks to the Markovian nature of the process $\a$, whence the solution $(V^\a,X^\a)$. Now by Feynman-Kac formula,  we see that $F$ satisfies the PDE:
\bea
	\label{eq60}
	F_t + \dfrac{1}{2}F_{vv}\sigma^2+\dfrac{1}{2}F_{xx} \rho^2+F_xb+F_x(\mu+ u \rho)= 0;\qq
		F(T,v,x)=v.
\eea
In light of the case $b=0$, we now 
look for the function $J(t,v,x)$ with the following properties:
\bea
	\label{eq38}
\left\{\ba{lll}
	J_x\rho(t,x) =  H(t,x) - F(t,v,x);\\
	 J_t+b(t,v,x)J_v+ \mu(t,x)J_x+\dfrac{1}{2}\sigma^2(t,v,x)J_{vv} + \dfrac{1}{2}\rho^2(t,x) J_{xx}=0,;\\
	J(T,v,x) \geq 0, \q\mbox{and}\q	J(T,v,x) = 0 \,\,{\mbox{\rm iff}}\,\, v = g(x).		
	\end{array}\right.
\end{eqnarray}
If such function $J$ exists, then a simple application of It\^o's formula shows that 
\beaa
	\hE^{\hP}[W^\a_T] = 	\hE^{\hP}[ -J(T,V_T,X_T) + J(0,v_0,x_0)] \leq J(0,v_0,x_0), 
\eeaa
and the equality holds when  $V_T = g(X_T)$ $\hP$-a.s., which would imply the optimality of the solution to STPBVP. To find such a function, we modify \eqref{eq10} as follows. Define
\begin{equation}
	\label{eq40}
	J(t,v,x) = \int_{g^{-1}(v)}^x \dfrac{H(t,y) - F(t,v,y)}{\rho(t,y)} dy + G(t,v):=\bar{J}(t,v,x)+G(t,v),
\end{equation}
where $G(t,v)$ is a function to be determined. 
We have the following result. 
	\begin{thm}
		\label{thm2}
		Assume Assumption \ref{assump1}, and assume further that  $b$ and $\si$ take the form:
$$b(t,v,x)=b_0(t)+b_1(t)v,\q \sigma(t,v,x)=\sigma(t,v).
$$
Then, a weak solution $(\bar\hP, \bar{V}, \bar{X}, \bar\a)$  to STPBVP (\ref{SDEVX})  with $\a$ having the affine structure  is an equilibrium strategy if 
 the following compatibility condition holds:
		\bea
		\label{eq61}
		\left\{\ba{lll}
			\rho(b_0(t)+b_1(t)H)+ \rho^2\Big[ (u_0+u_1v) \pa_xF  - \pa_xH (u_0+ u_1 H)\Big] = 0;\\
			\label{eq59}
			-\pa_x\mu\rho+\pa_x\rho\mu+\pa_t\rho+\dfrac{1}{2}\rho^2\pa_{xx}\rho = 0.
		\ea\right.
		\eea
where $H$ and $F$ satisfy the PDEs (\ref{eq58}) and (\ref{eq60}), respectively.
	\end{thm}
	
	{\it Proof.}
We shall argue that under the compatibility condition (\ref{eq61}), the function $J$ defined by (\ref{eq40}) satisfies the desired properties (\ref{eq38}). To see this, first note that the first identity in (\ref{eq38}) is trivial by definition of the function $J$ (\ref{eq40}). Next, we observe that $J(T,v,x) = \bar{J}(T,v,x)+G(T,v)=\int_{g^{-1}(v)}^x \dfrac{g(y)- v}{\rho(t,y)} dy + G(T,v)$, where $\bar{J}(T, v,x)\ge 0$ and $\bar{J}(T,v,x)=0$ if and only if $x=g^{-1}(v)$, as we argued in Theorem \ref{thm1}. Therefore the function $J$ defined by (\ref{eq40}) satisfies the terminal condition in (\ref{eq38}) provided $G(T,v)\equiv 0$. 

		Let us now look at the PDE in \eqref{eq38}.  Plugging \eqref{eq40} into the PDE in \eqref{eq38}, we have 
\bea
\label{eq42}
0&=&\pa_t G+b(t,v)\pa_v G+\dfrac{1}{2}\sigma^2(t,v)\pa_{vv}G+
			\pa_t \bar{J}+b(t,v)\pa_v\bar{J}+\dfrac{1}{2}\sigma^2(t,v)\pa_{vv}\bar{J}\nonumber\\
&&+ 	\mu\dfrac{H - F}{\rho}+ \dfrac{1}{2}[(\pa_xH-\pa_xF)\rho - (H- F)\pa_x \rho].
\eea
We shall now find conditions that $G(t,v)$ is a desired function by showing that  the 
PDE (\ref{eq42}) with the terminal condition  $G(T,v) \equiv 0$ is independent of $x$.
To see this we simply	
take the partial derivative with respect to $x$ on both sides of \eqref{eq42},   multiply by $\rho^2$, and set it to zero to get (suppressing
variables):
\bea
\label{eq46}
0&=&(H-F)\Big[ -\pa_t\rho+ \pa_x\mu \rho - \mu \pa_x\rho -\dfrac{1}{2}\rho^2\rho_{xx}\rho\Big]  \\
&&- \rho\Big\{ \pa_tF+b\pa_v F+\mu\pa_xF+\dfrac{1}{2}\sigma^2\pa_{vv}F+\dfrac{1}{2}\rho^2\pa_{xx}F\Big\}+ \rho(\pa_tH +\mu\pa_xH+\dfrac{1}{2}\rho^2\pa_{xx}H). \nonumber
\eea

Now, since $H$ and $F$ satisfies PDEs  \eqref{eq58} and \eqref{eq60}, respectively, we can easily check that a sufficient condition for
(\ref{eq46}) is the compatibility condition (\ref{eq61}), proving the theorem. 
\qed
	\begin{rem}
	{\rm	By \eqref{eq40}, it seems that the function $J$ depends on the choice of the strategy $\a$ since both PDEs (\ref{eq58}) and (\ref{eq60}) (for $H$ and $F$) do. However, the PDE in \eqref{eq38} for $J$, as well as its terminal condition are independent of $u$. Consequently, the function $J$ is independent of the choice of $u$. 
\qed}
	\end{rem}


\begin{thebibliography}{1} 

\bibitem{ABO}
Aase, K. K.,  Bjuland T., and Oksendal, B., {\it Strategic insider trading equilibrium: a filter theory approach}, Afr. Mat., 23 (2012), pp.145-162.

\bibitem{Aronson-67} Aronson, D.G., {\it  Bounds for the fundamental solution of a parabolic equation}, Bull. Amer. Math. Soc.,
73 (1967), pp.890-896.

\bibitem{Aronson-68} Aronson, D.G., {\it Non-negative solutions of linear parabolic equations}, Ann. Scuola Norm. Sup. Pisa., 22 (1968), pp.607-694.

\bibitem{KB92}
Back, K., {\it Insider trading in continuous time}, The Review of Financial Studies, 5 (1992), pp.387-409.

\bibitem{BP98}
Back, K., and Pedersen, H., {\it Long-lived information and intraday patterns}, Journal of financial markets, 1 (1998), pp.385-402.

%\bibitem{BC08}
%Bain, A., and Crisan, D. {\sl Fundamentals of Stochastic Filtering}. Vol. 60. Springer Science \& Business Media, 2008.

\bibitem{FB02}
Baudoin, F., {\it Conditioned stochastic differential equations: theory, examples and application to finance}, Stochastic Processes and their Applications, 100(1-2) (2002), pp.109-145.

\bibitem{BPSV}
Barrios, B., Peral, I., Soria, F., and Valdinoci, E., {\it A Widder's type theorem for the heat equation with non-local diffusion}, (2013), arXiv:1302.1786v1.

\bibitem{Ben}
Bensoussan, A., {\sl Stochastic control of partially observable systems}, Cambridge University Press, (2004).


\bibitem{BHMO12}
Biagini, F., Hu, Y., Meyer-Brandis, T., and Oksendal, B., {\it Insider trading equilibrium in a market with memory}, Math. Finan. Econ., 6 (2012), pp.229-247.

\bibitem{BLM} 
Buckdahn, R., Li, J.,  and Ma, J.,  {\it A Mean-field Stochastic Control Problem with   Partial  Observations}, The Annals of Applied Probab., 27(5) (2017), pp.3201-3245. 

\bibitem{CalSt}
Caldentey, R., and Stacchetti, E., {\it  Insider Trading with a Random Deadline}, Econometrica, 1 (2010), pp.245-283.

\bibitem{CCD11}
Campi, L., \c{C}etin, U., and Danilova, A., {\it Dynamic Markov bridges motivated by models of insider trading}, Stochastic Processes
and their Applications, 121 (2011), pp.534-567.

\bibitem{CCD13a}
Campi, L., \c{C}etin, U., and Danilova, A., {\it Explicit construction of a dynamic Bessel bridge of dimension 3}, Electron. J. Probab.,
18 (2013), pp.1-25.

\bibitem{CCD13b}
Campi, L., \c{C}etin, U., and Danilova, A., {\it Equilibrium model with default and dynamic insider information}, Finance and Stochastics, 17 (2013), pp.565-585.

\bibitem{CCD21}
\c{C}etin, U., and Danilova, A., {\it On pricing rules and optimal strategies in general Kyle-Back models
}, (2021), arXiv:1812.07529v2.	

\bibitem{CEL}
Cocquemas, F., Ekren, I., and Lioui, A., {\it A general solution method for insider problems}, (2020), arXiv:2006.09518v1.


\bibitem{Dani}
Danilova, A., {\it Stock market insider trading in continuous time with imperfect dynamic information}, Stochastics, 82 (2010), pp.111-131.

\bibitem{Fer}
Fernique, X., {\it Intégrabilité des vecteurs gaussiens}, CR Acad. Sci. Paris Serie A 270, (1970), pp.1698-1699.

\bibitem{FPY}
Fitzsimmons, P., Pitman, J., and Yor, M., {\it Markovian bridges: construction, Palm interpretation, and splicing}, In Seminar on Stochastic Processes, (1992), pp.101-134. 

\bibitem{FW}
Fleming, W.H., and Rishel, R.W., {\sl Deterministic and stochastic optimal control}, Springer Science \& Business Media, 1 (2012).

\bibitem{FI93}
F\"ollmer, H. and Imkeller, P., {\it Anticipation cancelled by a Girsanov transformation
: a paradox on Wiener space}, Annales de l'I. H. P., section B, 29(4) (1993), pp.569-586.

\bibitem{FWY}
F\"ollmer, H., Wu, C.T., and Yor, M., {\it Canonical decomposition of linear transformations of two independent Brownian motions motivated by models of insider trading}, Stochastic Processes and Their Applications, 84 (1999), pp.137-164.

\bibitem{FKK}
Fujisaki, M., Kallianpur, G., and Kunita, H., {\it Stochastic differential equations for the non linear filtering problem}, Osaka Journal of Mathematics, 9(1) (1972), pp.19-40.

\bibitem{G62}
Guenther, R.B., {\it An Extension of Widder's Theorem}, Master Thesis, Oregon State University, (1962).

\bibitem{HS}
Holden, C.W., and Subrahmanyam, A., {\it Long-lived private information and imperfect Competition}, Journal
of Finance, 1 (1992), pp.247-270.

%\bibitem{IW}
%Ikeda, N. and Watanabe, S.,  {\sl Stochastic differential equations and diffusion processes}. Elsevier (2014).

\bibitem{K85}
Kyle, A. S., {\it Continuous auctions and insider trading}, Econometrica: Journal of the Econometric Society, (1985), pp.1315-1335.

\bibitem{LSM}
Lepeltier, J.P., and San Martin, J., {\it Backward stochastic differential equations with continuous coefficient}, Statistics \& Probability Letters, 32(4) (1997), pp.425-430.

%\bibitem{Liu}
%Liu, Y.  {\sl  Some reinsurance/investment optimization problems for general risk reserve
%models}. Ph.D. thesis, Purdue Univ. (2004).

\bibitem{MaLiu}
Liu, Y., and Ma, J., {\it BSDE with jumps and multiplicative coefficients}, Preprint (2005).

\bibitem{MSZ}
Ma, J., Sun, R., and Zhou, Y., {\it Kyle--Back Equilibrium Models and Linear Conditional Mean-Field SDEs}, SIAM Journal on Control and Optimization, 56(2) (2018), pp.1154-1180.

\bibitem{MWZZ}
Ma, J., Wu, Z., Zhang, D., and Zhang, J., {\it On Wellposedness of Forward-Backward SDEs --- A Unified Approach}, Ann. Appl. Probab, 25(4) (2015), pp.2168-2214.

\bibitem{MYbook}
Ma, J., and Yong, J., {\sl Forward-backward stochastic differential equations and their applications}, Springer Science \& Business Media No.1702, (1999). 

\bibitem{PR}
Protter, P., {\sl Integration and Stochastic Differential Equations—A New Approach}, (1990).

\bibitem{TanY}
Tan, Y., {\sl  Generalized Kyle-Back Equilibrium Models with Dynamic Information and Related Topics}, PhD dissertation, University of Southern California, (2022).

\bibitem{W44}
Widder, D., {\it Positive temperatures on an infinite rod}, Transactions of the American Mathematical Society, 55 (1944), pp.85-95.

\bibitem{W53}
Widder, D., {\it Positive temperatures on a semi-infinite rod}, Transactions of the American Mathematical Society, 75 (1953), pp.510-525.

\bibitem{WU}
Wu, C.T., {\it Construction of Brownian motions in enlarged filtrations and their role in mathematical models of insider trading}, (1999).


\end{thebibliography}
\end{document}